\documentclass[12pt,twoside,a4paper]{article}
\usepackage[english]{babel}
\usepackage[utf8]{inputenc}
\usepackage{amsmath}
\usepackage{amsfonts}
\usepackage[pagebackref]{hyperref}
\usepackage{amsthm}
\usepackage{thmtools}
\usepackage{amssymb}
\usepackage{bm}
\usepackage{tikz-cd}
\usepackage{xr}
\usepackage[nameinlink]{cleveref}
\usepackage{url}

\usepackage{enumitem}
\usepackage{authblk}
\usepackage[title]{appendix}
\usepackage{chngcntr}
\usepackage{apptools}
\usepackage{mathtools}
\usepackage{appendix}

\usepackage{bookmark}
\usepackage{appendix}
\definecolor{crimsonglory}{rgb}{0.75, 0.0, 0.2}
\definecolor{darkpowderblue}{rgb}{0.0, 0.2, 0.6}
\hypersetup{
	colorlinks   = true, 
	urlcolor     = darkpowderblue, 
	linkcolor    = darkpowderblue, 
	citecolor   = crimsonglory 
}

\theoremstyle{plain}
\newtheorem{theorem}{Theorem}[section] 

\newtheorem{definition}[theorem]{Definition} 
\newtheorem{prop}[theorem]{Proposition}
\newtheorem{cor}[theorem]{Corollary}
\newtheorem{lemma}[theorem]{Lemma}

\newtheorem{remark}[theorem]{Remark}
\newtheorem{remarks}[theorem]{Remarks}
\newtheorem{conj}[theorem]{Conjecture}
\newtheorem{assumption}[theorem]{Assumption}

\DeclareMathOperator{\GL}{GL}
\DeclareMathOperator{\SL}{SL}

\DeclareMathOperator{\spec}{Spec}

\DeclareMathOperator{\aut}{Aut}

\DeclareMathOperator{\fin}{fin}

\newcommand{\Sum}[2]{\displaystyle\sum_{#1}^{#2}}

\newcommand{\Z}{\mathbb{Z}}
\newcommand{\N}{\mathbb{N}}
\newcommand{\C}{\mathbb{C}}
\newcommand{\Q}{\mathbb{Q}}

\newcommand{\A}{\mathbb{A}}

\newcommand{\G}{\mathbb{G}}

\newcommand{\CE}{\mathcal{E}}
\newcommand{\CX}{\mathcal{X}}
\newcommand{\ld}{,\ldots,}
\newcommand{\CP}{\mathcal{P}}
\newcommand{\CO}{\mathcal{O}}

\usepackage[OT2,T1]{fontenc}
\DeclareSymbolFont{cyrletters}{OT2}{wncyr}{m}{n}
\DeclareMathSymbol{\Sha}{\mathalpha}{cyrletters}{"58}

\usepackage{tocloft}

\makeatletter
\newcommand{\address}[1]{\gdef\@address{#1}}
\newcommand{\email}[1]{\gdef\@email{\url{#1}}}
\newcommand{\@endstuff}{\par\vspace{\baselineskip}\noindent\small
	\begin{tabular}{@{}l}\scshape\@address\\\textit{E-mail address:} \@email\end{tabular}}
\AtEndDocument{\@endstuff}
\makeatother

\title{Zilber-Pink in $Y(1)^n$: Beyond multiplicative degeneration}
\author{Georgios Papas}
\date{}
\address{Faculty of Mathematics and Computer Science\\
	The Weizmann Institute of Science\\
	234 Herzl Street,Rehovot 76100, Israel}
\email{georgios.papas@weizmann.ac.il}

\begin{document}
	\maketitle
	
	\begin{abstract}We establish Large Galois orbits conjectures for points of unlikely intersections of curves in $Y(1)^n$, upon assumptions on the intersection of such curves with the boundary $X(1)^n\backslash Y(1)^n$, in the Zilber-Pink setting. As a corollary, building on work of Habegger-Pila and Daw-Orr, we obtain new cases of the Zilber-Pink conjecture for curves in $Y(1)^n$.
 \end{abstract}
	
	
\section{Introduction}

The main objective of our exposition is to establish lower bounds for the size of Galois orbits of points in curves in the moduli space $Y(1)^n$ coming from unlikely intersections of our curves with special subvarieties of $Y(1)^n$. These results, known as ``Large Galois orbits conjectures'' in the general field of unlikely intersections, constitute the main difficulty in establishing the validity of unlikely intersections results using the Pila-Zannier method.

The main application of the results we obtain is some cases of the Zilber-Pink conjecture for curves in $Y(1)^n$. The general strategy to establish the Zilber-Pink conjecture in this setting is due to Habegger and Pila, see \cite{habeggerpila1}, where the authors reduce the conjecture to a Large Galois orbits conjecture. Their main unconditional result is the following:

\begin{theorem}[\cite{habeggerpila1}, Theorem $1$]\label{habpila}Let $C\subset Y(1)^n$ be an irreducible curve defined over $\bar{\Q}$ that is asymmetric and not contained in a proper special subvariety of $Y(1)^n$. 
	
	Then the Zilber-Pink conjecture holds for $C$.
\end{theorem}

In the process of establishing \Cref{habpila}, Habegger and Pila also reduce the conjecture for any curve $C$ as above, without the asymmetricity condition, to establishing finiteness of points of intersection of our curve with so called ``strongly special'' subvarieties of the moduli space $Y(1)^n$. These will be subvarieties that are defined by equations of the form $\Phi_M(x_{i_1},x_{i_2})=\Phi_N(x_{i_3},x_{i_4})=0$ where $1\leq i_j\leq n$ are such that the sets $\{i_1,i_2\}\neq \{i_3,i_4\}$ and $i_1\neq i_2$, $i_3\neq i_4$.

Using this circle of ideas, Daw and Orr establish the following:
\begin{theorem}[\cite{daworr4}, Theorem $1.3$]\label{daworrthm}Let $C\subset Y(1)^n$ be an irreducible curve defined over $\bar{\Q}$ that is not contained in a special subvariety of $Y(1)^n$ and is such that its compactification $\bar{C}$ in $X(1)^n$ intersects the point $(\infty\ld\infty)$.
	
	Then the Zilber-Pink conjecture holds for $C$.
\end{theorem}

Either of the conditions, i.e. the ``asymmetricity condition'' of Habegger-Pila or the condition about the type of the intersection of the curve with the boundary $X(1)^n\backslash Y(1)^n$, is needed in order to establish the aforementioned ``Large Galois orbits conjecture''. In \cite{habeggerpila1} this is achieved via a height bound due to Siegel and N\'eron, for which the asymmetricity condition is crucial. On the other hand, in \cite{daworr4}, Daw and Orr employ Andr\'e's G-functions method to arrive at certain height bounds at the points of interest. These in turn imply the lower bound on the size of the Galois orbits once coupled with the isogeny estimates of Masser-W\"ustholtz, see \cite{masserwu}. It is this same method introduced by Andr\'e that we use here to go beyond the condition of Daw and Orr about the intersections of our curve with the boundary $X(1)^n\backslash Y(1)^n$.

To state our main result in the direction of Zilber-Pink we first introduce a bit of notation.

Let $C\subset Y(1)^n$, where $n\geq 2$, be a smooth irreducible curve defined over $\bar{\Q}$ and let $\bar{C}$ be its Zariski closure in $X(1)^n$. Identifying the modular curve $Y(1)$ with $\A^1$ and $X(1)$ with $\mathbb{P}^1$ it is easy to see that the boundary $X(1)^n\backslash Y(1)^n$ of $Y(1)^n$ in its compactification $X(1)^n$ can be written as a disjoint union of subvarieties of the form $\infty^{h}\times (\A^1)^{n-h}$, up to permutation of the coordinates, with $1\leq h\leq n$.

Let now $s_0\in \bar{C}(\bar{\Q})\backslash Y(1)^n$ be a fixed point in the boundary $X(1)^n\backslash Y(1)^n$. 

\begin{definition}\label{defcoord}Let $C$, $s_0$ be as above and let $\pi_i:X(1)^n\rightarrow X(1)$ denote the coordinate projections. 
	
	The coordinate $i$ will be called \textbf{smooth for} $C$ if $\pi_i(s_0)\in Y(1)$. A smooth coordinate $i$ for the curve $C$ will be called a \textbf{CM coordinate} for $C$ if in addition $\pi_i(s_0)$ is a CM point in $Y(1)$. Finally, the coordinate $i$ will be called \textbf{ singular for }$C$ if it is not smooth, i.e. if $\pi_i(s_0)=\infty$.\end{definition}

\begin{theorem}\label{zpmain}
	Let $C\subset Y(1)^n$ be an irreducible curve defined over $\bar{\Q}$ that is not contained in any proper special subvariety of $Y(1)^n$. Assume that $C$ is such that all but at most one of its coordinates are singular and its one possibly smooth coordinate is CM. 
	
	Then the Zilber-Pink conjecture holds for $C$.
\end{theorem}

We note here that our condition on intersection with the boundary is clearly distinct from that of \Cref{daworrthm}. On the other hand, arguing similarly to the examples following Theorem $1.3$ of \cite{daworr4} one can see that this condition is distinct from the ``asymmetricity condition'' of \Cref{habpila}. In other words, our \Cref{zpmain} does not follow from \Cref{habpila}.

For our most general Zilber-Pink-type statement see \Cref{zpgen}. In the case where $n=3$ we can establish slightly stronger results. Namely, as a corollary of \Cref{zpgen} the following holds:
\begin{theorem}
	Let $C\subset Y(1)^3$ be an irreducible curve defined over $\bar{\Q}$ that is not contained in a proper special subvariety of $Y(1)^3$. Assume that the curve intersects the boundary $X(1)^3\backslash Y(1)^3$ in a point which up to permutation of coordinates is of the form $(\infty, \zeta_1,\zeta_2)$ or $(\infty,\infty,\zeta_1)$ with $\zeta_1$, $\zeta_2$ singular moduli.
	
	Then the Zilber-Pink conjecture holds for $C$.
\end{theorem}

For an up to date account on the Zilber-Pink conjecture we point the interested reader to \cite{pilabook}.
\subsection{Summary}

In \Cref{section:background} we start with summarizing some technical aspects of Andr\'e's G-functions method in the specific case of $Y(1)^n$. In particular, we discuss how to associate a family of G-functions to the geometric objects we study, i.e. to morphisms $f:X\rightarrow S$ with $S$ a curve and $X$ a semiabelian variety of a type relevant to the Zilber-Pink in this setting. To accomplish this we need to make several technical assumptions about our elliptic schemes and the base curves these are over. We treat these technicalities in detail in \Cref{section:heightreduction} where we show that they do not impede us in any way in establishing the height bounds we are need.

We continue in \Cref{section:isogrels}, where we construct relations among archimedean values of our family of G-functions at points $s\in S(\bar{\Q})$ over which the fiber of the morphism $f$ above reflects an intersection in the moduli space $Y(1)^n$ with strongly special subvarieties. In other words, the case pertinent to the open cases of the Zilber-Pink Conjecture.

After this, in \Cref{section:heightbounds} we start by establishing the height bounds needed to deduce our Large Galois orbits statements based on our previous work in \cite{papaseffbrsieg}. We finish our exposition in \Cref{section:applizp} by noting down the Large Galois orbits statement we need and recording some examples of Zilber-Pink type statements that follow readily from these coupled with the results of \cite{habeggerpila1} and \cite{daworr4}.\\

\textbf{Acknowledgments:} The author thanks Yves Andr\'e for answering some questions about his work on G-functions and Chris Daw for interesting discussions around Andr\'e's G-functions method. The author is grateful to Bijay Bhatta and Martin Orr for pointing out errors in the statements of \Cref{heightboundzp}, \Cref{lgozp}, and \Cref{zpgen} in an earlier version of this paper. Finally, the author is grateful to the anonymous referee for their many useful remarks, suggestions, and corrections on the material of the original manuscript.

Throughout this work the author was supported by Michael Temkin's ERC Consolidator Grant 770922 - BirNonArchGeom. During the revision process the author received funding by the European Union (ERC, SharpOS, 101087910), and by the ISRAEL SCIENCE FOUNDATION (grant No. 2067/23). 
\subsection{Notation}\label{section:notation}

Given a number field $L$ and $y:=\Sum{i=0}{\infty}a_ix^i\in L[[x]]$ a power series, we write \begin{itemize}
	\item $\Sigma_L$ (resp.  $\Sigma_{L,\infty}$ and $\Sigma_{L,f}$) for the set of places (resp. archimedean and finite places) of $L$.
	
	\item for $v\in \Sigma_L$, $\C_v$ will denote the corresponding complete, algebraically closed field and $\iota_v:L\hookrightarrow \C_v$ will denote the natural embedding,
	
	\item we write $\iota_v(y(x))$ for the power series $\Sum{i=0}{\infty}\iota_v(a_i)x^i\in\C_v[[x]]$.	
\end{itemize}

For a family $y_j\in L[[x]]$ and $\iota_v$ as above, we let $R_v(\{y_1\ld y_N\}):=\max\{ R_v(\iota_v(y_j))\}$, where $R_v(f)$ is the radius of convergence of $f\in \C_v[[x]]$.

Finally, for $r>0$ we write $\Delta_r$ for the open complex analytic disc of radius $r$ and centered at $0$.
	
	
	\section{Recollections on the G-functions method}\label{section:background}
	
	Here we summarize the main technical setting of the G-functions method in $Y(1)^n$. For more details see $\S 2$ of \cite{papaseffbrsieg}. In the process we also fix the setting and notation that will be used for the majority of our exposition.
	
	For an introduction to the theory of G-functions we point the interested reader to \cite{dwork} and the original work of Andr\'e \cite{andre1989g}.
	
	\subsection{Families of tuples of elliptic curves}\label{section:setting}
	Here we collect some terminology and general facts about N\'eron models that appear throughout this section. From now on unless explicitly stated all base curves are affine.\\
	
	Let us consider a smooth geometrically irreducible curve $S'$ over a number field $K$, $s_0\in S'(K)$ a fixed point, and set $S:=S'\backslash\{s_0\}$.
	We also consider an abelian scheme of the form $f:\CX=\CE_1\times_S \ldots\times_S \CE_n\rightarrow S$, defined over $K$ as well, with $f_i:\CE_i\rightarrow S$ defining an elliptic curve over $S$ for all $i$. 
	
	Given $\CE_i\rightarrow S$ as above, we may naturally consider the N\'eron model $N_S(\CE_i)\rightarrow S'$ of $\CE_i$ over $S'$. By this we mean that we consider the N\'eron model of the generic fiber $(\CE_{i})_{\eta}/K(S)$, which will exist over some smooth projective curve $\bar{S}$ containing $S'$ by \cite{neron} $\S 1.4$ Theorem $3$, and then restrict this to $S'$. 
	
	By Grothendieck's semi-stable reduction theorem, see \cite{neron} $\S 7.4$ Theorem $1$ for our main reference here, there exists some finite cover $C'$ of $S'$ over which $\CE_i$ acquires semi-stable reduction. In other words, the N\'eron model, which we denote by $N_{C}(\CE_i)$, of $(\CE_{i})_{\eta}\times_{K(S)}K(C')$ will be such that all of its fibers have connected components that are either elliptic curves or $1$-dimensional tori, see \cite{neron} $\S 7.4$ for our terminology here. 
	
	In this case, we may write $f_i':\CE'_i\rightarrow C'$ for the connected N\'eron model of  $(\CE_{i})_{\eta}\times_{K(S)}K(C')$ over $C'$, by this we mean the group subscheme of $N_C(\CE_i)$ whose fiber over every point is the connected component of the identity of the respective fiber of $N_C(\CE_i)$. Note that by the semi-stability assumption each of these fibers will be either a torus or an elliptic curve.
	
	Repeating the above process $n$ times, one for each $f_i$, we end up with a finite cover $C'$ of $S'$ and a morphism\begin{center}
		$f':\CX':=\CE'_1\times_S \ldots \times_S  \CE'_n\rightarrow C'$,
	\end{center}which will be the fiber product of the various $f'_i$ as above.
	
	With \Cref{defcoord} in mind we introduce the following:
	\begin{definition}\label{defsmorsingcoordrel}Let $S'$, $s_0$, and $f'$ be as above. The coordinate $i$ is said to be\begin{enumerate}
			\item smooth for $S'$ if $(f'_i)^{-1}(c_0)$ is an elliptic curve,
			\item CM for $S'$ if $(f'_i)^{-1}(c_0)$ is furthermore a CM elliptic curve, and
			\item singular for $S'$ if it is not smooth, so that $(f'_i)^{-1}(c_0)\simeq \G_m$,
		\end{enumerate}where $c_0$ ranges over the points of $C'$ that are over $s_0$. \end{definition}
	
	\begin{remark}\label{remarksingular}
		Let $i$ be a singular coordinate in our terminology and consider $f'_i:\CE_i'\rightarrow C'$ as above and let $\bar{f}_i:\bar{\CE}'_i\rightarrow C'$ denote the compactification of $f'_i$, see for example \cite{nagata}. The closed points $c_0$ in $C'$ that are over $s_0$ will then be singular values of the morphism $\bar{f}_i$.
		
		Suppose one were not, let $C:=C'\times_{S'}S$ and replace for simplicity $C'$ by $C\cup \{c_0\}$. Since $\CE'_i|_C=\bar{\CE}_i'|_C$ is smooth over $C$ we get that $\bar{\CE}_i'\rightarrow C'$ will be smooth proper, and since $C'$ itself is smooth, $\bar{\CE}'_i$ would be a regular proper model of $\CE_{i,\eta}\times_{K(S)}K(C')$ over $C'$. 
		
		Now let $f_i^{\min}:\CE^{\min}\rightarrow C'$ be a minimal regular proper model of $\CE_{i,\eta}\times_{K(S)}K(C')$ over $C'$, for the existence of this see for example \cite{silverad} Theorem $IV.4.5$. The above would imply, since $\bar{\CE}_i'$ and $\CE^{\min}$ are $C'$-isomorphic by minimality of $\CE^{\min}$, that $f^{\min}$ is also smooth. But then, by \cite{silverad} Theorem $IV.6.1$, we get that $N_C(\CE_i)=\CE^{\min}$. This would in turn imply that $N_C(\CE_i)\rightarrow C'$ is smooth proper and hence that $\CE_{i,\eta}\times_{K(S)}K(C')$ has good reduction on $c_0$, see for example Proposition $2$ in $\S 1.4$ of \cite{neron}. This clearly contradicts our assumption.
	\end{remark}

	\subsection{Height bounds}\label{section:heightreduction}
	
	Our primary objective, from a technical standpoint, is to establish height bounds on curves supporting tuples of elliptic curves. We phrase these here in the form of the following conjecture.
	
	\begin{conj}\label{conjhtbound}
		Let $C$ be a smooth irreducible curve defined over $\bar{\Q}$ and $f_k:\CE_k\rightarrow C$ be non-isotrivial elliptic curves, where $1\leq k\leq n$. Assume that the morphism $C\rightarrow Y(1)^n$ associated to this family is such that its image is not contained in a proper special subvariety. 
		
		Then there exist positive constants $c_1$ and $c_2$, depending on the $f_k$ and $C$, such that for all $s\in C(\bar{\Q})$ for which there exist $i_1\ld i_4\in \{1\ld n\}$ with $i_1\neq i_2$, $i_3\neq i_4$, and $\{i_1,i_2\}\neq \{i_3,i_4\}$ such that $\CE_{i_1,s}$ is isogenous to $\CE_{i_2,s}$ and $\CE_{i_3,s}$ is isogenous to $\CE_{i_4,s}$, we have \begin{equation}
			h(s)\leq c_1\cdot [\Q(s):\Q]^{c_2},
		\end{equation}where $h$ is some Weil height on $C$.
	\end{conj}
	
	\subsubsection{Reductions}
	
	It is cases of \Cref{conjhtbound} that we establish under conditions on the morphisms $f_k$. We do this following the G-functions method of Y. Andr\'e. In order to be able to apply this method to our problem we will need to make several assumptions on the curve $C$ and the morphisms $f_k$. 
	
	These assumptions turn out to be weak enough in the sense that establishing \Cref{conjhtbound} under these is equivalent to the full conjecture. We collect these in a series of lemmas here. 
	
	\begin{lemma}\label{reduction1}
		It suffices to establish \Cref{conjhtbound} under the following additional assumptions:\\
		
		there exists a smooth irreducible curve $S'$ over $\bar{\Q}$ and a point $s_0\in S'(\bar{\Q})$ for which $S=S'\backslash \{s_0\}$ and furthermore such that the following hold\begin{enumerate}	
			\item\label{l1it1} the N\'eron model, $N_S(\CE_i)$, of each $\CE_i$ over $S'$ has semi-stable reduction for all $i$,
			
			\item\label{l1it2} there exists at least one $j$ for which $N_S(\CE_j)$ has bad reduction at $s_0$,
			
			\item\label{l1it3} for all $i$ if the connected component of the identity of $N_S(\CE_i)_{s_0}$ is an elliptic curve, then it has everywhere semi-stable reduction, and 
			
			\item\label{l1it4} if $\iota:K\hookrightarrow \C$ is some complex embedding and $k$ is a singular coordinate we have that the local monodromy around $s_0$ acts unipotently on $R^1(f_k)_{*}\Q$.
		\end{enumerate}
	\end{lemma}
	\begin{proof} 
		For convenience let us fix one of the $1\leq k\leq n$ for now. 
		
		Let $\bar{S}$ be a smooth projective model over $\bar{\Q}$ of our $S$ and let $N(\CE_k)$ be the N\'eron model of $\CE_k$ over $\bar{S}$. Assume that $N(\CE_k)$ had potentially good reduction everywhere. Then we could find a smooth projective curve $\bar{C}$ that covers $\bar{S}$ such that the N\'eron model of the pullback, denoted $N_C(\CE_k)$, would be such that the connected component at each fiber over $\bar{C}$ is an elliptic curve. Letting $\CE_{k,C}'$ denote the associated connected N\'eron model, we would get therefore an elliptic curve $\CE_{k,C}'\rightarrow \bar{C}$. This would give a non-constant morphism $\bar{C}\rightarrow \mathbb{A}\hookrightarrow \mathbb{P}^1$ with $\infty$ not in its image, which is a contradiction. 
		
		Therefore we may find $s_0\in \bar{S}$ over which $N_S(\CE_k)$ does not have potentially good reduction. Note furthermore that we may assume that $s_0\in \bar{S}(\bar{\Q})$ since $\bar{S}$ and $S$ are both defined over $\bar{\Q}$ and $s_0\in \bar{S}\backslash S$. 
		
		Repeating the arguments of \Cref{section:setting}, we get a finite cover $\bar{C}$ of $\bar{S}$ such that, for $\eta$ the generic point of $\bar{S}$, all of the $\CE_{j,\eta}\times_{K(S)}K(\bar{C})$ have semi-stable reduction everywhere on $\bar{C}$. 
		
		Now, let us pick a point $c_0\in \bar{C}(\Q)$ over the point $s_0$ above so that without loss of generality $N_C(\CE_1)$ has bad reduction over $c_0$. Let us set $K:=\Q(\bar{C},c_0)$ to be a number field of definition of $\bar{C}$ and $c_0$. Again by Grothendieck's semi-stable reduction theorem, \cite{neron} Theorem $1$ of $\S 7.4$, there exists a finite extension $L_j/K$ such that if the connected component of the identity, denoted $E_{j,0}$ for simplicity, of the fiber over $c_0$ of $N_{C}(\CE_j)$ is an elliptic curve, then $E_{j,0}\otimes_KL_j$ has everywhere semi-stable reduction.
		
		Setting $L$ to be the compositum of all such $L_j$ and base changing by $L$ we end up with a smooth projective curve $\bar{C}_{\fin}:=\bar{C}\times_K L$, where $K$ here is some number field of definition of $\bar{C}$, and a dense open subset $C\subset \bar{C}_{\fin}$ such that the pair $(C,c_0)$ satisfy all of the conditions above.
		
		Assume now that \Cref{conjhtbound} holds for this pair. Note that the morphism $C_{\fin}\rightarrow \bar{C}\rightarrow \bar{S}$ is finite, with its degree depending only on the $f_k$ and $S$. Now note that if $s\in S(\bar{\Q})$ satisfies the hypothesis of the conjecture so will any of the points in $C_{\fin}$ that are over it. This is due to the fact that the elliptic curves on $C$ are obtained by base change of the original ones over $S$, after possibly removing finitely many points from $C$. The finiteness of the covering morphism together with classical properties of Weil's height machine, see for example \cite{hindrysilverman} Theorem $B.3.2$, imply that \Cref{conjhtbound} holds for $S$.\end{proof}
	
	Working under the assumptions of \Cref{reduction1}, as in \Cref{section:setting}, we denote by $f'_i:\CE'_i\rightarrow S'$ the connected N\'eron model of of ``$\CE_i$ over $S'$'', as defined there.
	\begin{lemma}\label{reduction2}
		It suffices to prove \Cref{conjhtbound} under the following additional assumptions:\\
		
		We have that $S=S'\backslash \{\xi_1\ld \xi_l\}$ where $S'$ is a smooth irreducible curve such that for each $1\leq t\leq l$ the pair $(S'_t:=S\cup \{\xi_t\}, \xi_t)$ satisfies the properties of \Cref{reduction1}. Let us also consider $K/\Q$ to be some finite extension over which $S'$ and all $\xi_t$ are all defined. We then, furthermore, assume that there exists a regular $\CO_K$-model $\mathfrak{S}$ of $S'$ as well as a semi-abelian scheme $\mathfrak{X}\rightarrow \mathfrak{S}$, and a rational function $x\in K(S')$ for which the following hold:

		\begin{enumerate}
			\item $\mathfrak{X}_K\simeq \CE_1'\times_S \ldots\times_S  \CE'_n$ as semi-abelian schemes over $S'$,
			
			\item $\{\xi_1\ld \xi_l\}=\{ s\in S'(\bar{\Q}); x(s)=0\}$ are simple zeroes of $x$, and 
			
			\item the morphism $x:S'\rightarrow \mathbb{P}^1$ induced from the above $x$ is Galois, or in other words the group $\aut_x(S'):=\{\sigma\in \aut(S'); x=x\circ \sigma\}$ acts transitively on the $\bar{\Q}$-fibers of $x$.\end{enumerate}
		
	\end{lemma}
	
	\begin{proof}Due to \Cref{reduction1}, we may assume that $S=S'\backslash \{s_0\}$ with $(S',s_0)$ that satisfies the assumptions there.
		
		Arguing as in Lemma $2.18$ of \cite{papaseffbrsieg} we may find a finite cover $\nu:C'\rightarrow S'$ where $C'$ satisfies all of the assumptions made here. We note that here the $\xi_t$ will be precisely the pre-images of $s_0$, this is due to the original construction of these covers by Daw and Orr, see \cite{daworr4} Lemma $5.1$. 
		
		Note again that by construction one has that the points $s\in S'(\bar{\Q})$ that are of interest in \Cref{conjhtbound} will have pre-images in $C'(\bar{\Q})$ that will also satisfy the same hypothesis. This is due to the fact that the corresponding semi-abelian schemes $(\CE'_i)_{C'}\rightarrow C'$ will just be the pullbacks of the $\CE_i'\rightarrow S'$ via the cover $\nu$. 
		
		Thus, we may conclude as in the proof of \Cref{reduction1} using basic properties of the Weil height machine.\end{proof}
	
	\subsection{G-functions and relative periods}\label{section:gfuns}
	
	From now on we assume that we are in the situation where the assumptions of both \Cref{reduction1} as well as \Cref{reduction2} are in place. In particular, as explained in the exposition of \Cref{section:setting}, we have semi-abelian schemes $f'_k:\CE_k'\rightarrow S'$ with $1$-dimensional fibers, where $1\leq k\leq n$ as per usual. We also fix from now on a number field $K$ over which everything is defined, i.e. the curve $S'$, the morphisms $f'_k$ as well as the points $\xi_t\in S'(\bar{\Q})$ in the notation of \Cref{reduction2}.
	
	As in \Cref{remarksingular} for those singular coordinates, we get that the aforementioned $\xi_t$ will be singular values of a certain compactification $\bar{f}_k$ of our morphism $f'_k$ as described there. For technical reasons needed in the proof of the following theorem we need to work under the following additional:
	\begin{assumption}\label{assummonodromy}Let $\iota:K\hookrightarrow \C$. For all singular coordinates $k$ we have that the local monodromy around any $s_0\in \{\xi_1\ld\xi_l\}$ acts unipotently on $R^1(f_k)_{\iota,*}\Q$.
	\end{assumption}
	
	\begin{remark}\label{rmkassumptionunipotent}Here $(f_{k})_{\iota}$ stands for the analytification of the morphism $f_k$ with respect to the embedding $\iota$.
		
		This additional assumption is harmless from the perspective of proving \Cref{conjhtbound}. Indeed, by the monodromy theorem, see for example \cite{landman}, the monodromy operator, acting on $R^1(f_k)_{\iota,*}\Q$, in a small archimedean disc centered at any of the above $s_0$ will act quasi-unipotently.
		
		We may thus find an \'etale cover of our $S'$, say $C'_{mon}$, such that once we base change the morphism $\bar{f}_k:\bar{\CE}_k\rightarrow S'$ with the \'etale cover, the local system over the analytification of $C'_{mon}$ corresponding to $R^1(f_k)_{\iota,*}\Q$ will have unipotent monodromy. Note here that by Riemann's existence theorem such an \'etale cover will be algebraic over $\bar{\Q}$. This might force us to work over some finite extension $L/K$.
		
		As in the proof of \Cref{reduction1}, again due to the properties of the Weil height, we can see that \Cref{assummonodromy} is indeed harmless in our setting. This is due to the fact that validity of \Cref{conjhtbound} for $S'$ is equivalent to the validity of the conjecture in any of its \'etale covers, or for its base change over $L/K$ as above.
		
		Finally, we note here that this \'etale base change, or the base change by a finite field extension $L/K$, does not alter the various assumptions about our setting made in either \Cref{reduction1} or \Cref{reduction2}.\end{remark}
	
	Let $x$ be as in \Cref{reduction2} and for simplicity let us fix $s_0$ to be one of the points in $\{\xi_\ld\xi_l\}$ described there. Then, under \Cref{assummonodromy}, and up to removing a finite set of points from the curve $S$, the following holds:
	\begin{theorem}[\cite{papaseffbrsieg} Theorem $2.7$]\label{gfunsthm}There exists a choice of a global basis $\{\omega_{2k-1},\omega_{2k}\}$ of $H^1_{DR}(\CE_k/S)$ and matrices of G-functions $Y_{G,k}=(y_{i,j,k})\in\GL_2 (\bar{\Q}[[x]])$ such that:\\
		
		for every archimedean embedding $\iota:\bar{\Q}\hookrightarrow \C$ and any simply connected $V\subset \Delta_{R_{\iota}(y_{i,j,k})}^{*}$, where the latter denotes a punctured disc contained in $S_{\iota}^{an}$ and centered at $s_0$, there exists a choice of local frames $\Gamma_{k,\iota}$ for the local systems $R^1(f_{k,\iota})_{*}\Q|_{V}$ with respect to which the relative period matrix $\CP$ of \begin{center}
			$H^1_{DR}(\CX/S)\otimes_{\C}\CO_{S^{an}_{\iota}}\rightarrow R^1(f_{\iota})_{*}\Q\otimes \CO_{S^{an}_\iota}$
		\end{center}over $V$ is block diagonal with blocks $\CP_k:=\begin{pmatrix}
			\frac{1}{2\pi i}\int_{\gamma_{k,j}}\omega_{2k-2+i}
		\end{pmatrix}_{1\leq i,j\leq 2}$, $1\leq k\leq n$, given by:
		\begin{enumerate}
			\item for all $k$ that are smooth for $S'$ and for all $s\in V$ we have $\CP_k(s)=\iota(Y_{G,k}(x(s))))\cdot \Pi_{k,\iota}$ for some $\Pi_k\in GL_2(\C)$, such that if, furthermore, $k$ is CM for $S'$ we have\begin{equation}
				\Pi_{k,\iota}=\begin{pmatrix}
					\frac{\varpi_{k,\iota}}{2\pi i}&0\\0&\varpi_{k,\iota}^{-1}.
				\end{pmatrix}
			\end{equation}
			
			\item for all $k$ that are singular for $S'$ and for all $s\in V$ we have \begin{equation}\label{eq:gfunssmooth}
				\CP_{k}(s)=	\iota(Y_{G,k}(x(s)))\cdot \begin{pmatrix}
					\iota(d_k)&e_{k,\iota}\\ \iota(d'_k)&e'_{k,i}
				\end{pmatrix} \cdot \begin{pmatrix}
					1&N_k\log \iota(x(s))\\0&1
				\end{pmatrix},
			\end{equation}where $d_k$, $d'_k\in \bar{\Q}$ are independent of the place $\iota$ and $N_k\in \Q$.
		\end{enumerate}
	\end{theorem}
	
	Throughout the text we will work under the following assumptions about our semiabelian scheme $f':\CX'\rightarrow S'$:
	\begin{assumption}\label{assumptionfield} The coefficients of the G-functions $y_{i,j,k}$ and the $d_k$, $d'_k$ are all in $K$.
	\end{assumption}
	
	In the case of singular coordinates we will also need the following:
	\begin{theorem}\label{lemmasymplectic}[\cite{papaseffbrsieg} Lemma $2.11$]For all $k$ singular coordinates for $S'$, there exists a choice of basis $\{\omega'_{2k-1}, \omega'_{2k}\}$ of $H^1_{DR}(\CE_k/S)$, possibly different from that in \Cref{gfunsthm}, such that with respect to the same frames $\Gamma_{k,\iota}$ the first column of the relative period matrix $\CP_{k,\iota}$, defined as in \Cref{gfunsthm}, is comprised of G-functions $y'_{i,1,k}$.\end{theorem}
	
	\subsection{Our setting}\label{section:goodcovers}
	
	Throughout the text we will work with G-functions as in \Cref{section:gfuns} associated to a point $\xi_t$ of our curve $S':=S\cup\{\xi_1\ld\xi_l\}$ which will satisfy the assumptions of \Cref{reduction1} as well as \Cref{reduction2}. The need to work with such a curve was first highlighted in \cite{daworr4}. 
	
	Unless otherwise stated we will also be working under the following ``Hodge genericity'' assumption about the abelian schemes $\CX:=\CX'|_S\rightarrow S$ that encodes our tuples of elliptic curves:
	\begin{assumption}\label{hodgegenericity}The image of the induced morphism $S\rightarrow Y(1)^n$ is a curve that is not contained in any proper special subvariety of $Y(1)^n$. 	 
	\end{assumption}

	The main property of this $S'$ is that it comes equipped with an $x\in K(S')$ which has simple zeroes only at the points $\xi_t$, $1\leq t\leq l$. It is important to note here that the ``Galois property'' of the cover $S'\rightarrow \mathbb{P}^1$ given by $x$ will guarantee that the notion of smooth, resp. CM, resp. singular, coordinates for $S'$ is independent of the point $\xi_t$ of reference.
	
	Let us write $S'_t:=S'\backslash \{\xi_j:j\neq t\}$, so that $S:=S'_t\backslash\{\xi_t\}$ is independent of the indexing $t$. Now the Galois property of the cover $x:S'\rightarrow \mathbb{P}^1$ guarantees that for each $1\leq t\leq l$ there exists $\sigma_t\in \aut_x(S')$ such that $\sigma_t(\xi_1)=\xi_t$. Pulling back the semi-abelian scheme $\CX'|_{S'_t}$ above via $\sigma_t$ we end up with a new semi-abelian scheme $\CX'_t\rightarrow S'_1$, one for each $1\leq t\leq n$.
	
	By \Cref{gfunsthm} and \Cref{lemmasymplectic} to each pair $(f'_t:\CX'_t\rightarrow S'_1,\xi_1)$ we get an associated family of G-functions \begin{center}
		$\mathcal{Y}_{t}:=\{y^{(t)}_{i,j,k}:k \text{ smooth for } S'\}\cup \{{y'}^{(t)}_{i,1,k}:k\text{ singular for }S'\}$.
	\end{center}
	
	The drawback of the above construction is that even though our original $\CX\rightarrow S$ was such that the induced morphism in $Y(1)^n$ had Hodge generic image this might no longer be true for the morphism induced by the abelian scheme $\CX_1\times_S\ldots\times_S\CX_l\rightarrow S$, where $\CX_t:=\CX'_t|_{S}$. 
	
	In order to work around this, following \cite{daworr4}, instead of working with all of the above G-functions we work with a subset of them. Namely, we consider $\Lambda\subset \{1\ld l\}$ that is a set of representatives for the equivalence relation \begin{center}
		$i\sim j$ if $(\CX_i)_{\eta}$ is isogenous to $(\CX_j)_{\eta}$
	\end{center}where $\eta$ stands for the generic point of $S_1$.
	
	To each $\lambda\in \Lambda$ we can then associate a single $\mathcal{Y}_{\lambda}$, say by associating the $\mathcal{Y}_t$ for which $t$ is minimal in the equivalence class $\lambda$. The family of G-functions we will be interested in applying Andr\'e-Bombieri's ``Hasse principle'' to, see Theorem $5.2$ in Chapter $VII$ of \cite{andre1989g}, will be the family $\mathcal{Y}=\{\mathcal{Y}_{\lambda} :\lambda\in \Lambda \}$.
	
	Under \Cref{hodgegenericity} we can then describe the so-called ``trivial relations'', see $VII. 5.1$ in \cite{andre1989g} for our terminology here, among our family $\mathcal{Y}$ of G-functions.
	\begin{theorem}[\cite{papaseffbrsieg}, Theorem $3.3$]\label{trivialrels}Assuming that \Cref{hodgegenericity} holds, we have that $\mathcal{Y}^{\bar{\Q}-Zar}$ is the subvariety of $\spec(\bar{\Q}[X^{(\lambda)}_{i,j,k}])$ defined by $Z(I_0)$ where \begin{equation}
			I_0:=\langle \det(X^{(\lambda)}_{i,j,k})-1:1\leq k\leq n, \text{ k is smooth for } S' \rangle,
		\end{equation}and where $X^{(\lambda)}_{i,j,k}$ are indeterminates in a polynomial ring $\bar{\Q}[X^{(\lambda)}_{i,j,k}]$ with the obvious restrictions on the indices, i.e. $\lambda\in\Lambda$, and the $i$, $j$, $k$ are such that $j=1$ if $k$ is a singular coordinate and $j=1$, $2$ otherwise.
	\end{theorem}
	
	Finally, we record a definition that we adopt throughout the exposition.
	\begin{definition}
		Any semiabelian scheme $f':\CX'\rightarrow S'$ as above, i.e. one that satisfies the assumptions of \Cref{reduction1,reduction2}, as well as \Cref{assummonodromy}, \Cref{assumptionfield}, and \Cref{hodgegenericity}, will be called $G_{ZP}$\textbf{-admissible}.
	\end{definition}

    
\section{Isogenies and archimedean relations}\label{section:isogrels}

We work in the general setting described in \Cref{section:goodcovers} which we consider fixed from now on. In particular, as we did in \Cref{section:goodcovers}, we assume throughout that \Cref{assummonodromy}, \Cref{assumptionfield}, and \Cref{hodgegenericity} hold for our curve $S'$.

For simplicity we will mostly work with a fixed $\lambda\in \Lambda$ in the setting of \Cref{section:goodcovers}. Therefore to ease our notation we will simply write $s_0$ for the point $\xi_1$ in the neighborhood of which we consider our G-functions and $y_{i,j,k}$ for the G-functions in question, instead of the more accurate $y^{(\lambda)}_{i,j,k}$. Furthermore, with notational simplicity in mind, from now on we also write ``$y_{i,1,k}$'' for the G-functions associated to singular coordinates as in \Cref{lemmasymplectic}, instead of ``$y_{i,1,k}'$''.

\subsection{Isogenies and periods} 

We record here the following lemma, which appears practically as Proposition $4.4$ of \cite{daworr4}. 
\begin{lemma}\label{isogenieslemma}Let $E_1$ and $E_2$ be elliptic curves defined over some number field $L$. Assume that there exists a cyclic isogeny $\phi:E_1\rightarrow E_2$ of degree $\deg(\phi)=M$ which is also defined over $L$.
	
	Let $P_k$ be the full period matrix of $E_k$, $k=1$, $2$, with respect to some fixed archimedean embedding $\iota:L\hookrightarrow \C$, some fixed bases $\{\gamma_{k,1},\gamma_{k,2}\}$ of $H_1(E_{k,\iota},\Z)$, and some fixed symplectic bases $\{\omega_{k,1},\omega_{k,2}\}$ of $H^{1}_{DR}(E_k/L)$ for which $\omega_{k,1}\in F^1:=e_k^{*}\Omega_{E_k/L}\subset H^1_{DR}(E_k/L)$ for $k=1$, $2$, where $e_k:L\rightarrow E_k$ stands for the respective zero section.
	
	Then, there exist $a$, $b$, $c\in L$ and $p$, $q$, $r$, $s\in\Z$ with $\det\begin{pmatrix}
		a&0\\b&c
	\end{pmatrix}=\det\begin{pmatrix}
		p&q\\r&s
	\end{pmatrix}=M$ such that \begin{equation}\label{eq:isogenieslemma}
		\begin{pmatrix}
			a&0\\b&c
		\end{pmatrix}\cdot P_1=P_2\cdot \begin{pmatrix}
			p&q\\r&s
		\end{pmatrix}
	\end{equation}
\end{lemma}
\begin{proof}Let $\omega_1\in F^1\backslash\{0\}$ and $\omega_2\in H^1_{DR}(E_1/L)$ so that set $\{\omega_1. \omega_2\}$ is a symplectic basis with respect to the polarizing form. Similarly let $\{\omega'_1,\omega'_2\}$ be a basis of $H^1_{DR}(E_2/L)$ with the same properties. Let also $\{\gamma_1,\gamma_2\}$ and $\{\gamma'_1,\gamma'_2\}$ be symplectic bases of $H_1(E_1,\Z)$ and $H_1(E_2,\Z)$ respectively. 
	
	We then have that there exists $a\in L$ such that $\phi^{*}(\omega'_1)=a\cdot \omega_1$ and there exist $b$, $c\in L$ such that $\phi^{*}(\omega'_2)=b\cdot \omega_1+c\cdot \omega_2$. On the one hand for the homology we know that there exist $p$, $q$, $r$, and $s\in \Z$ such that $\phi_{*}(\gamma_1)=p\cdot \gamma'_1+r\cdot \gamma'_2$ and $\phi_{*}(\gamma_2)=q\cdot \gamma'_1+s\cdot \gamma'_2$. 
	
	On the other hand we have \begin{equation}
		\int_{\gamma_j}^{} \phi^{*}(\omega'_i)= \int_{\phi_{*}(\gamma_j)}^{}\omega'_i.
	\end{equation}
	
	Combining this with the above we obtain for $i=j=1$ \begin{center}
		$a\int_{\gamma_1}^{}\omega_1= p \cdot \int_{\gamma'_1}^{}\omega'_1+r \cdot \int_{\gamma'_2}^{}\omega'_1$,
	\end{center}and similar relations from the other pairs of indices. Their combination is just the above equality of matrices.\end{proof}
\subsection{The toy case: $n=3$} \label{section:toycase}

Zilber-Pink for curves starts taking meaning for $n\geq 3$. In this subsection we work with the minimal such dimension, i.e. here $n=3$.\\

We write $\CX'_0$ for the connected fiber at $s_0$ of the N\'eron model of $\CX$ over $S'$, as per our usual notation explained in \Cref{section:setting}. Keeping in mind that we are working under the assumptions of \Cref{reduction1} and \Cref{reduction2}, we get that, up to some permutations of the coordinates, there are in practice three cases that we need to examine in this scenario:
\begin{enumerate}
	\item $\CX'_0=\G_m^3$ this has been dealt with in \cite{daworr4},
	
	\item $\CX'_0=\G_m^2\times E$ with $E$ some elliptic curve, or
	
	\item $\CX'_0=\G_m \times E\times E'$ with $E$ and $E'$ (not necessarily distinct) elliptic curves. 
\end{enumerate}

It is special cases of cases $2$ and $3$ above that we are interested in. In what follows we shall keep notation as above for the decomposition of the fiber $\CX_0$. Namely we shall assume, which we can do without loss of generality, that the potentially singular coordinates for $S'$ are the first two. We refer to each of the cases by the type of fiber that appears over $s_0$. 

Throughout this subsection we fix a point $s\in S(\bar{\Q})$. We write $E_1\times E_2\times E_3$ for the fiber $\CX_{s}$ at $s$ of our family and assume that there exist $\phi_1:E_3\rightarrow E_1$ and  $\phi_2:E_3\rightarrow E_2$ cyclic isogenies of degree $\deg(\phi_k)=M_k$. We also let $L_s$ be the compositum of $K(s)$ with the fields of definition of these isogenies.

Let us assume that $s$ is $v$-adically close to $\xi_t$ with respect to some fixed  archimedean place $v\in\Sigma_{L_s,\infty}$, i.e. that $s$ is in the connected component of $x^{-1}(\Delta_{\min\{1,R_{\iota_v}(y_{i,j,k})\}})$ that contains $\xi_t$. Leaving the general argument for the proof of \Cref{archrelastionsisogenies}, we assume for now that $t=\lambda$ is a minimal element in its equivalence class, for the equivalence relation defined in \Cref{section:goodcovers}.

In such a case, we have that the fiber $(\CX_{t})_{\sigma_t^{-1}(s)}$, of the abelian scheme $\CX_t$ introduced in \Cref{section:goodcovers}, will also be of the same type. It is with these considerations in mind that we work from now. In particular we look to apply \Cref{isogenieslemma} to construct relations among the values at $x(\sigma_t^{-1}(s))=x(s)$ of the family $\mathcal{Y}_\lambda$. As noted earlier in this simplified version where our $\lambda$ is fixed we drop any mention of it.

\begin{definition}\label{atypicisogenous}
	$1.$ Any point $s\in S (\bar{\Q})$ as above will be called a point with \textbf{unlikely isogenies} for the semiabelian scheme $f':\CX'\rightarrow S'$.\\
	
	$2.$ We call the field $L_s$ defined above the\textbf{ field of coefficients of the point} $s$.
\end{definition}

By \Cref{gfunsthm} we have three matrices of G-functions, one for each coordinate, where as mentioned earlier we will, for ease of notation, write $Y_{G,k}(x)$ rather than the more accurate ``$Y^{(\lambda)}_{G,k}(x)$'' to denote them. For convenience we also write $Y_{G,k}(x(s))= \begin{pmatrix}
	\tilde{h}^{(k)}_{i,j}
\end{pmatrix}$ for the entries of these matrices, i.e. the values of the G-functions at $\xi:=x(s)$.

We will also write $\CP_{k}(s)$ for the values at $s$ of the respective relative period matrices of $f'_{t,k}:\CE'_{t,k}\rightarrow S'_t$ constructed with respect to the bases and trivializations used in \Cref{section:goodcovers} to construct the family $\mathcal{Y}$ associated to $s_0$.

\subsubsection{The case $\G_m^2\times E$}\label{section:gm2xe}

From \Cref{isogenieslemma} we get that there exist $a_k,$ $b_k$, $c_k\in L$ and $p_k$, $q_k$, $r_k$ $s_k\in \Z$ such that \begin{equation}\label{eq:gm2crosscm0}
	\begin{pmatrix}a_k&0\\b_k&c_k\end{pmatrix}\Pi_3 \begin{pmatrix}\tilde{h}^{(3)}_{i,j}	\end{pmatrix}\begin{pmatrix}\varpi_{0,1}&\varpi_{0,2}\\ \varpi_{0,3}&\varpi_{0,4}	\end{pmatrix}=\Pi_k\begin{pmatrix} \tilde{h}^{(k)}_{i,j}	\end{pmatrix}   \begin{pmatrix} d_k&e_k\\d'_k&e'_k	\end{pmatrix}\begin{pmatrix} 1&N_k\log(\xi)\\0&1	\end{pmatrix}\begin{pmatrix}p_k&q_k\\r_k&s_k	\end{pmatrix}.
\end{equation}Here $\Pi_3$ is the change of basis matrix from the basis $\{\omega_{5,s},\omega_{6,s}\}$ of $H^1_{DR}(\CE_{3,s}/L)$ constructed in \Cref{gfunsthm} to the basis used in \Cref{isogenieslemma} and $\Pi_k:=\Pi_{k,1}\cdot \Pi_{k,2}$, for $k=1$, $2$, is the product of the change of basis matrices $\Pi_{k,2}$, that passes from the basis of $H^1_{DR}(\CE_{k,s}/L)$ chosen in \Cref{gfunsthm} to that given by \Cref{lemmasymplectic}, and $\Pi_{k,1}$, which passes from the basis of $H^1_{DR}(\CE_{k,s}/L)$ chosen in \Cref{lemmasymplectic} to that chosen in \Cref{isogenieslemma}.

Note here that $d_k$, $d_k'\in K$ by \Cref{assumptionfield}. To ease our notation a little we set $e_{0,k}:=d_kN_k\log(\xi)+e_k$ and $e'_{0,k}:=d'_kN_k\log(\xi)+e'_k$. Also, writing  $\Pi_k\begin{pmatrix}\tilde{h}^{(k)}_{i,j}\end{pmatrix}=\begin{pmatrix}{h}^{(k)}_{i,j}\end{pmatrix}$, we may rewrite the above in the more useful form \begin{equation}\label{eq:gm2crosscm1}
	\begin{pmatrix}a_k&0\\b_k&c_k\end{pmatrix} \begin{pmatrix}h^{(3)}_{i,j}	\end{pmatrix}\begin{pmatrix}\varpi_{0,1}&\varpi_{0,2}\\ \varpi_{0,3}&\varpi_{0,4}	\end{pmatrix}=\begin{pmatrix} h^{(k)}_{i,j}	\end{pmatrix}   \begin{pmatrix} d_k&e_{0,k}\\d'_k&e'_{0,k}	\end{pmatrix} \begin{pmatrix}p_k&q_k\\r_k&s_k	\end{pmatrix}.\end{equation}

\begin{remark}[The CM case]As we will see, the case where $E$ is CM is easier to handle. The vanishing of the periods $\varpi_{0,2}$ and $\varpi_{0,3}$ turns out to make computations of relations feasible.
	
	We record here for our convenience \eqref{eq:gm2crosscm1} under the assumption that $E$ is CM:
	\begin{equation}\label{eq:gm2crosscm2}
		\begin{pmatrix}a_k&0\\b_k&c_k\end{pmatrix} \begin{pmatrix}{h}^{(3)}_{i,j}	\end{pmatrix}\begin{pmatrix} \frac{\varpi_{0,3}}{2\pi i}&0\\0& \varpi_{0,3}^{-1}	\end{pmatrix}=\begin{pmatrix}{h}^{(k)}_{i,j}	\end{pmatrix} \begin{pmatrix} d_k&e_{0,k}\\d'_k&e'_{0,k}	\end{pmatrix}  \begin{pmatrix}p_k&q_k\\r_k&s_k	\end{pmatrix}
	\end{equation}
\end{remark}

\subsubsection*{Towards relations}

There are two potential ways to go from \eqref{eq:gm2crosscm1} to relations among the $h^{(k)}_{i,j}$. They both use the same technique inspired from \cite{daworr4} Proposition $4.4$. The first of these will only use the G-functions $y^{(t)}_{i,1,k}$ corresponding to the first column of the matrices $\begin{pmatrix} h^{(k)}_{i,j}	\end{pmatrix}\begin{pmatrix} d_k&e_{0,k}\\d'_k&e'_{0,k}	\end{pmatrix}$ coming from the two singular coordinates.

Here we have chosen to work in the greatest possible generality for two reasons. First of all, these computations appear throughout all cases we will deal with in one way or another. Secondly, the computations themselves reveal the limitations of current methods at least to the knowledge of the author.\\

\textbf{First way:} Multiply both sides of \eqref{eq:gm2crosscm1} on the left by the vector \begin{equation}\label{eq:intermediatvecs1}
	(g_1^{(k)}, g_2^{(k)}):=(d_kh^{(k)}_{2,1}+d'_kh^{(k)}_{2,2}, -(d_kh^{(k)}_{1,1}+d'_kh^{(k)}_{1,2})),
\end{equation}to get  the following 
\begin{equation}\label{eq:gm2generalway11}
	(a_kg^{(k)}_1+b_kg^{(k)}_{2},c_kg^{(k)}_2)\begin{pmatrix}h^{(3)}_{i,j}	\end{pmatrix}\begin{pmatrix}\varpi_{0,1}&\varpi_{0,2}\\ \varpi_{0,3}&\varpi_{0,4}	\end{pmatrix}=(0,-\frac{D_k}{2\pi i})\begin{pmatrix}p_k&q_k\\r_k&s_k	\end{pmatrix},
\end{equation}where $D_k:=\det(\Pi_k)\in L_s^{\times}$.

Here we are using the fact that $\det\begin{pmatrix} \tilde{h}^{(k)}_{i,j}	\end{pmatrix}\begin{pmatrix} d_k&e_k\\d'_k&e'_k	\end{pmatrix}\begin{pmatrix} 1&N_k\log(\xi)\\0&1	\end{pmatrix}=\det\begin{pmatrix} d_k&e_{0,k}\\d'_k&e'_{0,k}	\end{pmatrix}=\frac{1}{2\pi i}$, from the Legendre relation, see \cite{langelliptic} Chapter $18$ $\S 1$ bottom of page $241$. On the other hand $\det\begin{pmatrix}\tilde{h}^{(k)}_{i,j}\end{pmatrix}=1$ for all $k$. To see this, note that by the Legendre relations the determinant of the relative period matrix is constant and equal to $\frac{1}{2\pi i}$. On the other hand, the matrices of G-functions will be such that $Y(0)=I_2$. This holds by construction, see for example the proof of Claim $3.7$ in \cite{papaszp} and the proof of Theorem $3.1$ there.

Setting \begin{equation}\label{eq:intermediatedvecs2}(H^{(k)}_1,H^{(k)}_2)= -D_k^{-1}((a_kg^{(k)}_1+b_kg^{(k)}_{2})h^{(3)}_{1,1}+ c_k g^{(k)}_2h^{(3)}_{2,1}, (a_kg^{(k)}_1+b_kg^{(k)}_{2})h^{(3)}_{1,2}+c_kg^{(k)}_2h^{(3)}_{2,2}), 
\end{equation}one gets that \begin{equation}\label{eq:gm2generalway12}(H_1^{(k)}, H^{(k)}_2)\cdot \begin{pmatrix}\varpi_{0,1}&\varpi_{0,2}\\ \varpi_{0,3}&\varpi_{0,4}	\end{pmatrix}=(\frac{r_k}{2\pi i},\frac{s_k}{2\pi i}).
\end{equation}

This finally translates to the pair of relations \begin{equation}\label{eq:generalway1}H^{(k)}_1 \varpi_{0,1}+H^{(k)}_2 \varpi_{0,3}=\frac{ r_k}{2\pi i} \text{ and }
	H^{(k)}_1 \varpi_{0,2}+H^{(k)}_2 \varpi_{0,4}=\frac{s_k}{2\pi i}
\end{equation}

\begin{remark}
	Note that the transcendence degree of the (possibly transcendental) periods $\varpi_{0,i}$ and $\pi$ over $\bar{\Q}$ is $\leq 4$ and conjecturally under Grothendieck's period conjecture will be equal to $4$ when our elliptic curve is not CM. In spirit we do not have enough equations to ``get rid of'' all of them and create a relation among the values of the $h^{(k)}_{i,j}$ with coefficients in $\bar{\Q}$.
\end{remark}

\textbf{Second way:} Here we are using all of the G-functions from the singular coordinates.\\

Multiply both sides of \eqref{eq:gm2crosscm1} on the left by the vector \begin{equation}\label{eq:intermediatvecs3}
	(h^{(k)}_{2,1}, -h^{(k)}_{1,1}),
\end{equation}using the fact that $\det\begin{pmatrix}h^{(k)}_{i,j}\end{pmatrix}=D_k$ for $k=1$, $2$, to get  the following 
\begin{equation}\label{eq:gm2generalway21}
	(a_kh^{(k)}_{2,1}-b_kh^{(k)}_{1,1},-c_kh^{(k)}_{1,1})\begin{pmatrix}h^{(3)}_{i,j}	\end{pmatrix}\begin{pmatrix}\varpi_{0,1}&\varpi_{0,2}\\ \varpi_{0,3}&\varpi_{0,4}	\end{pmatrix}=(0,-D_k)   \begin{pmatrix} d_k&e_{0,k}\\d'_k&e'_{0,k}	\end{pmatrix}\begin{pmatrix}p_k&q_k\\r_k&s_k	\end{pmatrix}.
\end{equation}

Setting \begin{equation}\label{eq:intermediatvecs3.5}
	(g_1^{(k)}, g_2^{(k)}):=(-D_k^{-1}(a_kh^{(k)}_{2,1}-b_kh^{(k)}_{1,1}), D_k^{-1}c_k h^{(k)}_{1,1}),\text{ and then}
\end{equation} 
\begin{equation}\label{eq:intermediatedvecs4}(H^{(k)}_1,H^{(k)}_2)= (g_1^{(k)}h^{(3)}_{1,1}  +g_2^{(k)} h^{(3)}_{2,1},g_1^{(k)}h^{(3)}_{1,2}+g_2^{(k)}h^{(3)}_{2,2}), 
\end{equation}one gets that \begin{equation}\label{eq:gm2generalway22}(H_1^{(k)}, H^{(k)}_2)\cdot \begin{pmatrix}\varpi_{0,1}&\varpi_{0,2}\\ \varpi_{0,3}&\varpi_{0,4}	\end{pmatrix}=(d'_k,e'_{0,k})\begin{pmatrix}p_k&q_k\\r_k&s_k	\end{pmatrix}.
\end{equation}

This finally translates to the pair of relations \begin{equation}\label{eq:generalway2a}H^{(k)}_1 \varpi_{0,1}+H^{(k)}_2 \varpi_{0,3}=d'_kp_k+e'_{0,k}r_k, \text{ and }
	H^{(k)}_1 \varpi_{0,2}+H^{(k)}_2 \varpi_{0,4}=d'_kq_k+e'_{0,k}s_k.
\end{equation}

Now repeat the above from the start by multiplying both sides of \eqref{eq:gm2crosscm1} on the left by the vector \begin{equation}
	(h^{(k)}_{2,2}, -h^{(k)}_{1,2}),
\end{equation}to get \begin{equation}\label{eq:gm2generalway23}
	(a_kh^{(k)}_{2,2}-b_kh^{(k)}_{1,2},-c_kh^{(k)}_{1,2})\begin{pmatrix}h^{(3)}_{i,j}	\end{pmatrix}\begin{pmatrix}\varpi_{0,1}&\varpi_{0,2}\\ \varpi_{0,3}&\varpi_{0,4}	\end{pmatrix}=(D_k,0)   \begin{pmatrix} d_k&e_{0,k}\\d'_k&e'_{0,k}	\end{pmatrix}\begin{pmatrix}p_k&q_k\\r_k&s_k	\end{pmatrix}.
\end{equation}
Setting
\begin{equation}\label{eq:intermediatvecs5}
	(g_3^{(k)}, g_4^{(k)}):=(D_k^{-1}(a_kh^{(k)}_{2,2}-b_kh^{(k)}_{1,2}),-D_k^{-1}c_kh^{(k)}_{1,2}),\text{ and then}
\end{equation} \begin{equation}\label{eq:intermediatvecs6}
	(H^{(k)}_3, H^{(k)}_4):= (h^{(3)}_{1,1}g^{(k)}_3+h^{(3)}_{2,1}g^{(k)}_4,h^{(3)}_{1,2}g^{(k)}_3+h^{(3)}_{2,2}g^{(k)}_4)
\end{equation}
one then keeps going as earlier to reach an analogue of \eqref{eq:generalway2a}, namely one gets:\begin{equation}\label{eq:generalway2b}
	(H^{(k)}_3 \varpi_{0,1}+H^{(k)}_4 \varpi_{0,3} ,H^{(k)}_3 \varpi_{0,2}+H^{(k)}_4 \varpi_{0,4})=    ( d_kp_k+e_{0,k}q_k,  d_kr_k+e_{0,k}s_k).
\end{equation}

\begin{remark}
	The advantage to the previous computations is evident. We now have more potential relations to try to create some relation strictly among the $H_i^{(k)}$ by eliminating the $\varpi_{0,j}$. The drawback is that through this way we have introduced more transcendental numbers, namely the $e_{0,k}$ and $e'_{0,k}$. 
	
	Nevertheless, this still seems to not be enough, at least to the author, to deal with the problem of creating archimedean relations among the $h^{(k)}_{i,j}$ unless we make assumption about the transcendental numbers that appear above.
\end{remark}

\subsubsection*{The subcase where $E$ has CM}

From now on assume that $E$, the fiber in the third coordinate of the fiber $\CX_0$, has CM, so we can use \eqref{eq:gm2crosscm2} instead. Using the same exact argument as the one employed in the ``First way'' of the previous paragraph, we get from \eqref{eq:generalway1} in this setting that

\begin{equation}\label{eq:gm2xcm1}H^{(k)}_1\frac{\varpi_{0,3}}{2\pi i}=\frac{r_k}{2\pi i}, \text{ and }
	H^{(k)}_2 \varpi_{0,3}^{-1}=\frac{s_k}{2\pi i}.
\end{equation}Multiplying these together we get $H^{(k)}_1\cdot H^{(k)}_2=\frac{r_ks_k}{2\pi i}$ for $k=1$, $2$. 

From this, one gets that either $H^{(k)}_j=0$ for some $k$ and $j$ or, if all of the $r_k$, $s_k\neq0$, that $H^{(1)}_1\cdot H^{(1)}_2r_2s_2=H^{(2)}_1\cdot H^{(2)}_2r_1s_1$. We thus have:
\begin{lemma}\label{localfactorse2cm}
	Let $f':\CX'\rightarrow S'$ be a $G_{ZP}$-admissible semiabelian scheme. Assume that\footnote{Due to the Galois property of the cover $x:S'\rightarrow \mathbb{P}^1$ if one of the fibers has this property so will the rest.} the fibers $\CX'_{\xi_t}$ are of $\G_m^2\times E$-type with $E$ CM.
	
	Let $s\in S(\bar{\Q})$ be some point with unlikely isogenies and let $L_s$ be its associated field of coefficients. Then if $s$ is $v$-adically close to $\xi_\lambda$ with respect to some archimedean place $v\in \Sigma_{L_s,\infty}$, there exists $R_{s,v}\in L_s[X^{(\lambda)}_{i,j,k}]$ such that the following hold\begin{enumerate}
		\item $\iota_v(R_{s,v}(\mathcal{Y}_{\lambda}(x(s))))=0$,
		
		\item $R_{s,v}$ is homogeneous of degree $\deg(R_{s,v})\leq 4$, and
		
		\item $R_{s,v}\notin I_0\leq L_s[X^{(\lambda)}_{i,j,k}]$, where $I_0$ is the ideal defined in \Cref{trivialrels}.
	\end{enumerate}
\end{lemma}
\begin{proof}From the above discussion we have that either $H^{(k)}_j=0$ for some $k$ and $j$, or that $H^{(1)}_1\cdot H^{(1)}_2r_2s_2=H^{(2)}_1\cdot H^{(2)}_2r_1s_1$.
	
	We start with some remarks. Note that by the discussion preceding \eqref{eq:gm2crosscm1} we have that by definition the first column of the matrix $\Pi_{k,2}\cdot (\tilde{h}_{i,j}^{(k)})\cdot   \begin{pmatrix} d_k&e_{0,k}\\d'_k&e'_{0,k}	\end{pmatrix}$ is nothing but\footnote{We remind the reader here that the $y_{i,1,k}$ with $k$ singular are those given in \Cref{lemmasymplectic}.} $\begin{pmatrix}
		\iota_v(y_{1,1,k}^{(\lambda)}(x(s)))\\  \iota_v(y_{2,1,k}^{(\lambda)}(x(s))).
	\end{pmatrix}$. Writing $\Pi_{k,1}=(a_{i,j,k})\in \SL_2(L_s)$ we thus have that the intermediate vector \eqref{eq:intermediatvecs1} is nothing but 
	
	\begin{equation}\label{eq:auxilgm2cm1}\begin{split}
			(g_1^{(k)}, g_2^{(k)})=  (\iota_v(a_{2,1,k}y_{1,1,k}^{(\lambda)}(x(s))        +a_{2,2,k} y_{2,1,k}^{(\lambda)}(x(s))    ) ,\\ -\iota_v(a_{1,1,k}y_{1,1,k}^{(\lambda)}(x(s))        +a_{1,2,k} y_{2,1,k}^{(\lambda)}(x(s)))).\end{split}\end{equation}
	
	Writing $\Pi_3:=(a_{i,j,3})$ we thus get that $h_{i,j}^{(3)}$ are linear combinations of the entries of the matrix $(\iota_v(y^{(\lambda)}_{i,j,3}(x(s))))$, which are by construction the values of G-functions we are interested in.
	
	Therefore the equations $H^{(1)}_j=0$ and $H^{(1)}_1\cdot H^{(1)}_2r_2s_2=H^{(2)}_1\cdot H^{(2)}_2r_1s_1$, will correspond to a polynomial $R_{s,v}$ that by construction will satisfy all but the final conclusion of our lemma. The rest of this proof focuses on this final part of our statement, i.e. the non-triviality of the $R_{s,v}$. 
	
	For notational convenience, we drop from now on any reference to $\lambda$, i.e. the index referring to the root $\xi_\lambda$ of the ``local parameter'' $x$ on our curve.\\
	
	\textbf{Case $1$: }$H^{(k)}_j=0$\\
	
	Let us assume without loss of generality that $H^{(1)}_1=0$, i.e. that $j=k=1$. Then $R_{s,v}$ will be the following polynomial\begin{equation}\label{eq:localfactorgm2cm1}\begin{split}
			R_{s,v}=-c_1(a_{2,1,3}X_{1,1,3}+a_{2,2,3}X_{2,1,3})(a_{1,1,1}X_{1,1,1}+a_{1,2,1}X_{2,1,1})\\
			+((a_1a_{2,1,1} -b_1  a_{1,1,1} )     X_{1,1,1}+(a_1a_{2,2,1}-b_1a_{1,2,1})X_{2,1,1})\cdot\\
			\cdot (   a_{1,1,3}X_{1,1,3}+a_{1,2,3}X_{2,1,3}  ).\end{split}	
	\end{equation}
	
	We first argue that at least one of the coefficients of the presentation of $R_{s,v}$ as a sum of monomials is non-zero. Assume that this is not so. Then looking at the coefficients of the monomials $X_{1,1,1}X_{1,1,3}$ and $X_{1,1,1}X_{2,1,3}$ we get that \begin{equation}\label{eq:coeff111}
		(a_1a_{2,1,1}-b_1a_{1,1,1})a_{1,1,3}-c_1a_{1,1,1}a_{2,1,3}=0\text{, and}
	\end{equation}
	\begin{equation}\label{eq:coeff211}
		(a_1a_{2,1,1}-b_1a_{1,1,1})a_{1,2,3}-c_1a_{1,1,1}a_{2,2,3}=0.
	\end{equation}

Since $\det(\Pi_3)\neq0$, the above implies that $(a_1a_{2,1,1}-b_1a_{1,1,1},-c_1a_{1,1,1})=(0,0)$. Note that $c_1\neq 0$ by construction thus $a_{1,1,1}=0$. This in turn gives $a_1a_{2,1,1}=0$ and since again $a_1\neq0$ we get $a_{2,1,1}=0$ which would imply $\det(\Pi_{1,1})=0$.

Let us now assume that $R_{s,v}\in I_0$. Working modulo the ideal \begin{center}
	$\langle X_{i,i,k}^{(\lambda')}-1, X_{j,j',k}^{(\lambda')}; \lambda'\neq \lambda, j\neq j'\rangle$,
\end{center}thus ignoring the issue of multiple $\lambda$'s from now on, we conclude $R_{s,v}\in I_0'$ where $I_0':=(\det(X_{i,j,3})-1 )\leq \bar{\Q}[X_{i,j,k};1\leq k\leq 3]$. 

Now note that $I_0'\subset \mathfrak{m}:=\langle X_{i,i,3}-1,X_{2,1,3}, X_{1,2,3}1\leq i\leq 2\rangle$. Writing $\phi:\bar{\Q}[X_{i,j,k};1\leq k\leq 3] \rightarrow \bar{\Q}[X_{i,j,k};1\leq k\leq 3] /\mathfrak{m}\simeq \bar{\Q}[X_{i,j,k};1\leq k\leq 2]$ for the quotient by $\mathfrak{m}$ our assumption implies that $\phi(R_{s,v})=0$. This would again imply \eqref{eq:coeff111}, that now appears as the coefficient of $X_{1,1,1}$ in $\phi(R_{s,v})$, as well as \eqref{eq:coeff211}, now the coefficient of $X_{2,1,1}$ in $\phi(R_{s,v})$, would both be zero. As we just saw this is impossible.\\

\textbf{Case $2$: }$H^{(1)}_1\cdot H^{(1)}_2r_2s_2=H^{(2)}_1\cdot H^{(2)}_2r_1s_1$\\

In this case we will have that $r_k$, $s_k\neq 0$ for $k=1$, $2$ by construction. Let us write $R_{H_j,k}$ for the polynomial corresponding to $H^{(k)}_j$, for example $-D_1\cdot R_{H_1,1}$ is the polynomial described in \eqref{eq:localfactorgm2cm1}. 

Then $R_{s,v}= r_2s_2R_{H_1,1}R_{H_2,1}-r_1s_1R_{H_1,2}R_{H_2,2}$. The same computations giving \eqref{eq:localfactorgm2cm1} give\begin{equation}\label{eq:localfactorgm2cm2}\begin{split}
		-D_1\cdot R_{H_2,1}=-c_1(a_{2,1,3}X_{1,2,3}+a_{2,2,3}X_{2,2,3})(a_{1,1,1}X_{1,1,1}+a_{1,2,1}X_{2,1,1})\\
		+((a_1a_{2,1,1}-b_1  a_{1,1,1} ) X_{1,1,1}     +(a_1a_{2,2,1}-    +b_1a_{1,2,1})X_{2,1,1})\cdot\\
		\cdot (   a_{1,1,3}X_{1,2,3}+a_{1,2,3}X_{2,2,3}  ).\end{split}	
\end{equation}

Writing \begin{center}
	$R_{H_1,1}=C_1 X_{1,1,1}X_{1,1,3}+C_2 X_{1,1,1} X_{2,1,3}+ C_3X_{2,1,1}X_{1,1,3}+C_4 X_{2,1,1}X_{2,1,3}$ 
\end{center}we notice that \begin{center}
	$R_{H_2,1}=C_1 X_{1,1,1}X_{1,2,3}+C_2 X_{1,1,1} X_{2,2,3}+ C_3X_{2,1,1}X_{1,2,3}+C_4 X_{2,1,1}X_{2,2,3}$, 
\end{center}i.e. the coefficients are the same with at least one of them being non-zero.

By symmetry one has\begin{center}
	$R_{H_1,2}=C'_1 X_{1,1,2}X_{1,1,3}+C'_2 X_{1,1,2} X_{2,1,3}+ C'_3X_{2,1,2}X_{1,1,3}+C'_4 X_{2,1,2}X_{2,1,3}$ 
\end{center}we notice that \begin{center}
	$R_{H_2,2}=C'_1 X_{1,1,2}X_{1,2,3}+C'_2 X_{1,1,2} X_{2,2,3}+ C'_3X_{2,1,2}X_{1,2,3}+C'_4 X_{2,1,2}X_{2,2,3}$, 
\end{center}i.e. the coefficients are again the same and at least one of them is non-zero.

Now, if $R_{s,v}\in I_0$, ignoring the rest of the $\lambda$'s by arguing as above, we would have $R_{s,v}\in \mathfrak{m}_1:=\langle X_{1,1,3}-1,X_{2,1,3},X_{1,2,3},X_{2,2,3}-1\rangle$. This in turn, by working modulo $\mathfrak{m}_1$ as earlier, implies that \begin{equation}\label{eq:contrag2cm}\begin{split}
		r_2s_2( C_1 X_{1,1,1}+C_3X_{2,1,1}) (C_2X_{1,1,1}+C_4X_{2,1,1})\\-r_1s_1(C'_1X_{1,1,2}+C_3'X_{2,1,2})(C_2'X_{1,1,2}+C_4'X_{2,1,2})=0.\end{split}
\end{equation}

The proof in the previous case shows that at least one of the $C_1$ and $C_2$, and similarly at least one of $C_3$ and $C_4$ are non-zero, and the same for the coefficients $C'_j$. If \eqref{eq:contrag2cm} were to hold we must then have that, without loss of generality, $C_2=C_4=0$. 

Then, noting that $I_0\subset \mathfrak{m}_2:= \langle X_{1,1,3}-1,X_{2,1,3},X_{1,2,3}-1,X_{2,2,3}-1\rangle$, we get  $R_{s,v}\in \mathfrak{m}_2$ which, again going back to $\bar{\Q}[X_{i,j,k};1\leq k\leq 3]$ and modding out by $\mathfrak{m}_2$, implies
\begin{equation}\label{eq:contrag2cm2}\begin{split}
		r_2s_2( C_1 X_{1,1,1}+C_3X_{2,1,1}) (C_1X_{1,1,1}+C_3X_{2,1,1})-F(X_{1,1,2},X_{2,1,2})=0,\end{split},
\end{equation}for some homogeneous polynomial $F\in \bar{\Q}[X_{1,1,2}, X_{2,1,2}]$. This is clearly impossible since $r_2s_2C_1\neq 0$ and the coefficient of $X_{1,1,1}^2$ is $r_2s_2C_1^2\neq 0$.\end{proof}
\subsubsection{The $\G_m\times E\times E'$ case}\label{section:gmee}

The same issue as in \Cref{section:gm2xe} pops up. Namely, there are too many possibly transcendental numbers that appear in our equations. Nevertheless, there are special cases here where we can extract relations among the values of the G-functions of our family.

\subsubsection*{$E'$ is CM}

Let us write $\begin{pmatrix}	\varpi_{0,1}&\varpi_{0,2}\\\varpi_{0,3}& \varpi_{0,4}\end{pmatrix}$ for the periods of the elliptic curve $E$ and $\begin{pmatrix}
	\frac{\varpi}{2\pi i}&0\\0&\varpi^{-1}\end{pmatrix}$ for those of $E'$. 

Working with the isogenous pair $\phi_2^{\vee}:E_2\rightarrow E_3$ of the fiber at $s$ we get the following, here as before we write $\Pi_k\cdot Y_{G,k}(x(s))=(h_{i,j}^{k})$, note that now $\Pi_k$ for $k=2$, $3$, are defined in the same manner as $\Pi_3$ in \Cref{section:gm2xe}:
\begin{equation}\label{eq:gmexcm1}
	\begin{pmatrix}a_3&0\\b_3&c_3\end{pmatrix} \begin{pmatrix}h^{(2)}_{i,j}	\end{pmatrix}\begin{pmatrix} \varpi_{0,i}	\end{pmatrix}=\begin{pmatrix} h^{(3)}_{i,j}	\end{pmatrix} \begin{pmatrix}\frac{\varpi}{2\pi i}&0\\0&\varpi^{-1}\end{pmatrix}  \begin{pmatrix}p_3&q_3\\r_3&s_3	\end{pmatrix}.\end{equation}

From this, one gets working as in the ``second way'' above \begin{equation}\label{eq:gmxexcm23a}
	(H^{(3)}_1\varpi_{0,1}+H^{(3)}_2\varpi_{0,3}, H^{(3)}_{1}\varpi_{0,2}+H^{(3)}_2\varpi_{0,4})=(\frac{r_3}{\varpi}, \frac{s_3}{\varpi}), \text{ and }
\end{equation}\begin{equation}\label{eq:gmxexcm23b}
	(H^{(3)}_3\varpi_{0,1}+H^{(3)}_4\varpi_{0,3}, H^{(3)}_{3}\varpi_{0,2}+H^{(3)}_4\varpi_{0,4})=(\frac{p_3\varpi}{2\pi i}, \frac{q_3\varpi}{2\pi i}).
\end{equation}

Now we look at the pair of isogenous elliptic curves $\phi:E_3\rightarrow E_1$. From the previous discussion, working as in the ``second way'' outlined in the previous section, we get: \begin{equation}
	\begin{pmatrix}a_1&0\\b_1&c_1\end{pmatrix}	\begin{pmatrix} h^{(3)}_{i,j}	\end{pmatrix} \begin{pmatrix}\frac{\varpi}{2\pi i}&0\\0&\varpi^{-1}\end{pmatrix}=\begin{pmatrix} h^{(1)}_{i,j}	\end{pmatrix} \begin{pmatrix} d_1&e_{0,1}\\d'_1&e'_{0,1}	\end{pmatrix}  \begin{pmatrix}p_1&q_1\\r_1&s_1\end{pmatrix}.\end{equation}

This will lead us to equations of the form $H^{(1)}_1\varpi =r_1$ and $H^{(1)}_2\frac{1}{\varpi}=\frac{s_1}{2\pi i}$. For the other pair of functions we get equations of the form \begin{equation}\label{eq:gmxexcm13b}
	(\frac{\varpi}{2\pi i}H^{(1)}_3, \frac{1}{\varpi}H^{(1)}_4)=(d_1p_1+e_{0,1}r_1, d_1q_2+e_{0,1}s_1).
\end{equation}

\begin{remark}These seem to not be sufficient for our purposes in dealing with the general case here, i.e. that where the other elliptic curve $E$ is generic. Once again, there are too many possibly transcendental numbers that appear in these equations.
\end{remark}

\subsubsection*{$E$ is also CM}

Let us now assume that $E$ is also CM. Then we can get relations in two different ways.\\

\textbf{First way:} Working as in \Cref{section:gm2xe}, namely the constructions under the assumption that $E$ is CM  there, we get from working with the isogenous pair $(E_1,E_3)$ the relations 

\begin{equation}
	H^{(1)}_1H^{(1)}_2=\frac{r_1s_1}{2\pi i}, 
\end{equation}and working with the pair $(E_1,E_2)$ we get the relation
\begin{equation}
	H^{(2)}_1H^{(2)}_2=\frac{r_2s_2}{2\pi i}.
\end{equation}From these we can get rid of $\pi$ and get a relation as before.\\

\textbf{Second way:} The second way is to work only with the pair $(E_2,E_3)$. One then gets a simplified version of the equation in \eqref{eq:gmexcm1}. Namely, one has:
\begin{equation}\label{eq:isogcmxcm1}
	\begin{pmatrix}a_3&0\\b_3&c_3\end{pmatrix} \begin{pmatrix}h^{(2)}_{i,j}	\end{pmatrix} \begin{pmatrix} \frac{\varpi'}{2\pi i}&0\\0&\varpi'^{-1}	\end{pmatrix}=\begin{pmatrix} h^{(3)}_{i,j}	\end{pmatrix} \begin{pmatrix}\frac{\varpi}{2\pi i}&0\\0&\varpi^{-1}\end{pmatrix}  \begin{pmatrix}p_3&q_3\\r_3&s_3	\end{pmatrix},\end{equation}where $\begin{pmatrix} \frac{\varpi'}{2\pi i}&0\\0&\varpi'^{-1}	\end{pmatrix}$ is the period matrix of $E$. 

We work much as in the ``second way'' outlined in \Cref{section:gm2xe}. Multiplying both sides of the above on the left by $(h^{(3)}_{2,2},-h^{(3)}_{1,2})$ we get:
\begin{equation}\label{eq:isogcmcmgeneralway21}
	(a_3h^{(3)}_{2,2}-b_3h^{(3)}_{1,2},-c_3h^{(3)}_{1,2})\begin{pmatrix}h^{(2)}_{i,j}\end{pmatrix}   \begin{pmatrix} \frac{\varpi'}{2\pi i}&0\\0&\varpi'^{-1}	\end{pmatrix}=(1,0)   \begin{pmatrix} \frac{\varpi}{2\pi i}&0\\0&\varpi^{-1}	\end{pmatrix}\begin{pmatrix}p_3&q_3\\r_3&s_3	\end{pmatrix}.
\end{equation}
As usual setting $(g_3,g_4):=(a_3h^{(3)}_{2,2}-b_3h^{(3)}_{1,2},-c_3h^{(3)}_{1,2})$ and  $(H_3,H_4):=(g_3h^{(2)}_{1,1}+g_4h^{(2)}_{2,1}, g_3h^{(2)}_{1,2}+g_4h^{(2)}_{2,2})$, we get \begin{equation}\label{eq:cmcmway2a}
	(\frac{\varpi'H_3}{2\pi i}, \frac{H_4}{\varpi'})=(\frac{\varpi p_3}{2\pi i},\frac{\varpi q_3}{2\pi i}).
\end{equation}

Multiplying \eqref{eq:isogcmxcm1} on the left on both sides by $(-h^{(3)}_{2,1},h^{(3)}_{1,1})$ and repeating the notation from earlier we end up with the relations:
\begin{equation}\label{eq:cmcmway2b}
	(\frac{\varpi'H_1}{2\pi i},\frac{H_2}{\varpi'})=(\frac{r_3}{\varpi},\frac{s_3}{\varpi}).\end{equation}

Combining this with \eqref{eq:cmcmway2a} gives \begin{equation}\label{eq:cmcmway2f}H_1H_2H_3H_4=p_3q_3r_3s_3.\end{equation}

\begin{remark}
	The $H_i$ correspond to homogeneous degree $2$ polynomials among the $h^{(k)}_{i,j}$. To turn \eqref{eq:cmcmway2f} into a relation coming from a homogeneous polynomial we can just multiply its right hand side by $1=\det(y_{i,j,3}(x(s)))^2\det(y_{i,j,2}(x(s)))^2$. 
	
	In what follows, these relations coming out of \eqref{eq:cmcmway2f} will be more useful. In other words, we will be working with the relations corresponding to the following equation \begin{equation}\label{eq:cmcmway2final}H_1H_2H_3H_4=\iota_v(p_3q_3r_3s_3\det(y_{i,j,3}(x(s)))^2\det(y_{i,j,2}(x(s)))^2).\end{equation}
\end{remark}

\begin{lemma}\label{localfactors2ndway}Let $f':\CX'\rightarrow S'$ be a $G_{ZP}$-admissible semiabelian scheme. Assume that $\CX'_{s_0}$, for any $s_0\in \{\xi_1\ld \xi_l\}$, is of $\G_m\times E\times E'$-type with $E$ and $E'$ CM.
	
	Let $s\in S(\bar{\Q})$ be some point with unlikely isogenies and let $L_s$ be its associated field of coefficients. Then if $s$ is $v$-adically close to $\xi_\lambda$ with respect to some archimedean place $v\in \Sigma_{L_s,\infty}$, there exists $R_{s,v}\in L_s[X^{(\lambda)}_{i,j,k}]$ such that the following hold\begin{enumerate}
		\item $\iota_v(R_{s,v}(\mathcal{Y}_{\lambda}(x(s))))=0$,
		
		\item $R_{s,v}$ is homogeneous of degree $\deg(R_{s,v})\leq 8$, and
		
		\item $R_{s,v}\notin I_0\leq L_s[X^{(\lambda')}_{i,j,k}:1\leq i,j\leq 2, k=2,3]$, where $I_0$ is the ideal defined in \Cref{trivialrels}.
	\end{enumerate}
\end{lemma}
\begin{proof}
	The first two properties follow by construction. The construction gives us two possible cases, either one of the entries of $\begin{pmatrix}
		p_3&q_3\\r_3&s_3
	\end{pmatrix}$ is zero, in which case we get a polynomial from the relation $H_j=0$, or they are all non-zero in which case we look to \eqref{eq:cmcmway2final} for the polynomial $R_{s,v}$ we want. 
	
	As in earlier proofs we are reduced to showing non-triviality of the relations in question, i.e. that $R_{s,v}\notin I_0$. As in previous proofs we drop any reference of the index ``$\lambda$'' from now on for notational simplicity. Note here that by modding out by the same ideal as in the proof of \Cref{localfactorse2cm} we may reduce to working in the polynomial ring $\bar{\Q}[X^{\lambda}_{i,j,k};1\leq i,j\leq 2, 1\leq k\leq 3]$, i.e. where we have gotten rid of the indeterminates with $\lambda'\neq\lambda$.\\
	
	\textbf{Case $1$: }$H_j=0$\\
	
	Without loss of generality assume that $H_3=0$, i.e. that $p_3=0$. Then by construction we have \begin{equation}
		R_{s,v}= C_1 X_{1,2,3}X_{1,1,2}+C_2 X_{1,2,3} X_{2,1,2} +C_3 X_{2,2,3}X_{1,1,2}+C_4 X_{2,2,3}X_{2,1,2},\text{ where}
	\end{equation}
	\begin{center}
		$C_1= (a_3a_{2,1,3}-b_3a_{1,1,3})a_{1,1,2}-c_3a_{1,1,3}a_{2,1,2}$, 
		
		$C_2= (a_3a_{2,1,3}-b_3a_{1,1,3})a_{1,2,2}-c_3a_{1,1,3}a_{2,2,2}$, 
		
		$C_3= (a_3a_{2,2,3}-b_3a_{1,2,3})a_{1,1,2}-c_3a_{1,2,3}a_{2,1,2}$, and 
		
		$C_4= (a_3a_{2,2,3}-b_3a_{1,2,3})a_{1,2,2}-c_3a_{1,2,3}a_{2,2,2}$. 
	\end{center}
	
	Now assume that $R_{s,v}=0$, i.e. that $C_j=0$ for all $j$. Since $\Pi_2=(a_{i,j,2})$ is invertible and by definition $c_3\neq 0$, we will then have that $a_{1,1,3}=a_{1,2,3}=0$ which clearly contradicts the fact that $\Pi_3=(a_{i,j,3})$ is invertible. Therefore, we get $R_{s,v}\neq0$ and the monomials in $R_{s,v}$ do not appear in the presentations of the two generators $\det(X_{i,j,2})-1$ and $\det(X_{i,j,3})-1$ of the ideal $I_0$, thus $R_{s,v}\notin I_0$ in this case.
	
	We note that furthermore, as in the proof of \Cref{localfactorse2cm}, we can see that at least one of $C_1$ and $C_2$ has to be non-zero and likewise for the pair $C_3$, $C_4$. Indeed, if say $C_1=C_2=0$, since $\det(\Pi_2)\neq 0$ and $a_3$, $c_3\neq 0$ we must have that $a_{1,1,3}=a_{2,1,3}=0$ but this once again contradicts the fact that $\det(\Pi_3)\neq0$.\\
	
	\textbf{Case $2$: }$H_1H_2H_3H_4=\iota_v(p_3q_3r_3s_3\det(y_{i,j,3}(x(s)))^2\det(y_{i,j,2}(x(s)))^2)$\\
	
	Again here we will have that all entries of the matrix $\begin{pmatrix}
		p_3&q_3\\r_3&s_3
	\end{pmatrix}$ are non-zero. Assume from now on that $R_{s,v}\in I_0$ and write $f_{k}:=\det(X_{i,j,k})-1$ for its two generators.
	
	Let us write $R_i$ for the polynomial corresponding to each of the $H_i$. In this sense we will have \begin{center}
		$R_{s,v}=R_1R_2R_3R_4-p_3q_3r_3s_3\det(X_{i,j,3})^4$.
	\end{center}

We have already seen that \begin{center}
	$R_3= C_1 X_{1,2,3}X_{1,1,2}+C_2 X_{1,2,3} X_{2,1,2} +C_3 X_{2,2,3}X_{1,1,2}+C_4 X_{2,2,3}X_{2,1,2}$.
\end{center}
Computing $R_4$ we see, as in the previous case, that we may write \begin{center}
	$R_4= C_1X_{1,2,3}X_{1,2,2}+C_2X_{1,2,3}X_{2,2,2}+C_3X_{2,2,3}X_{1,2,2}+C_4X_{2,2,3}X_{2,2,2}$,
\end{center}where $C_j$ are the exact same coefficients as above. 

Similar computations give \begin{center}
	$R_1=C_1' X_{1,1,3}X_{1,1,2} +   C_2'  X_{1,1,3}X_{2,1,2}     +C'_3    X_{2,1,3}X_{1,1,2}   +C_4' X_{2,1,3}X_{2,1,2}$ \text{ and }
\end{center}
\begin{center}
	$R_2=C_1'  X_{1,1,3} X_{1,2,2}+   C_2' X_{1,1,3}X_{2,2,2}     +C'_3    X_{2,1,3}X_{1,2,2}   +C_4' X_{2,1,3}X_{2,2,2}$,
\end{center} again with the same coefficients.

Let us first consider the ideal \begin{center}
	$\mathfrak{m}_1:=\langle f_2, X_{1,2,3},X_{2,1,3},X_{1,1,3}X_{2,2,3}-1 \rangle $,
\end{center}noting that $I_0\subset \mathfrak{m}_1$. Then from $R_{s,v}\in I_0$, working modulo $\mathfrak{m}_1$ as in the proof of \Cref{localfactorse2cm}, we conclude that the polynomial \begin{equation}\begin{split}
		Q_1:=(C_1'X_{1,1,2}+C_2'X_{2,1,2})  (C_1' X_{1,2,2}+C_2'X_{2,2,2})\\ (C_3X_{1,1,2}+C_4X_{2,1,2})(C_{3}X_{1,2,2}+C_4 X_{2,2,2})\end{split}
\end{equation}is such that $Q_1\in (f_2)$,where $(f_2)$ here denotes the principal ideal of $L_s(X_{1,1,3})[X_{i,j,2}:1\leq i,j\leq2]$. 

Noting that $(f_2)\subset \mathfrak{m}_2:=\langle X_{2,2,2}-1,X_{1,1,2}-1, X_{1,2,2},X_{2,1,2}\rangle$, we can see, looking at $Q_1$ modulo $\mathfrak{m}_2$, that \begin{equation}
	C_1'C_2'C_3C_4=0.
\end{equation}Let us assume without loss of generality that $C_1'=0$. Then from the discussion in the first part, from the symmetry of the definition of the $H_j$, we know that $C_2'\neq0$.

On the other hand, looking at $Q_1$ modulo any of the ideals $\mathfrak{m}_{3,n}:=\langle X_{2,2,2}-1,X_{1,1,2}-1, X_{1,2,2},X_{2,1,2} -n\rangle$, for which $I_0\subset \mathfrak{m}_{3,n}$ for all $n\in \N$, we see that \begin{equation}
	(C_2')^2(C_3+nC_4)C_4=0
\end{equation}holds for all $n\in \N$. This clearly implies that $C_4=0$ and hence $C_3\neq 0$ by our remarks in the previous case of the proof. 

Now the relations $C_1'=C_4=0$ imply \begin{center}
	$(C_2'C_3)^{-2}\cdot Q_1= X_{2,1,2}X_{2,2,2}X_{1,1,2}X_{1,2,2}\in (f_2)$,
\end{center} the latter viewed as an ideal in the ring $L_s(X_{1,1,3})[X_{i,j,2}:1\leq i,j\leq2]$. Since $(f_2)$ is prime this would imply that $X_{i,j,2}\in (f_2)$ for some pair $i$, $j$ which is clearly absurd.\end{proof}

\begin{remark}
	The distinct advantage of \Cref{localfactors2ndway} is that one only needs one isogeny to create the relations in question! The negligible for our arguments disadvantage is that one has that the polynomial will be of higher degree potentially than the one that would correspond to relations coming out of the ``first way'' described above.
\end{remark}
\subsection{Archimedean relations at points with unlikely isogenies}

Let $s\in S(L)$ and let $v\in \Sigma_{L,\infty}$ where $S$ as usual is some curve satisfying the assumptions in \Cref{reduction1,reduction2}. Then we will say that $s$ is ``$v$-adically close to $0$'' if $|x(s)|_v<\min\{1, R_{\iota_v}(y^{(\lambda)}_{i,j,k})\}$.

Putting everything together from the previous section we can conclude with the following:
\begin{prop}\label{archrelastionsisogenies}
	
	Let $f':\CX'\rightarrow S'$ be a $G_{ZP}$-admissible semiabelian scheme of relative dimension $n$. Let $s\in S(\bar{\Q})$ be a point that has unlikely isogenies and assume that not all of the isogenous coordinates of $\CX_s$ are singular for $S'$, while all of these that are smooth coordinates for $S'$ are furthermore CM for $S'$.
	
	Then, there exists a homogeneous polynomial $R_{s,\infty}\in L_s[X_{i,j,k}: 1\leq i,j,\leq 2, 1\leq k\leq n]$, where $L_s/K(s)$ is a finite extension, such that the following hold: \begin{enumerate}
		
		\item $\iota_v(R_{s,\infty}(\mathcal{Y}(x(s))))=0$ for all $v\in \Sigma_{L_s,\infty}$ for which $s$ is $v$-adically close to $0$,
		
		\item $[L_s: \Q]\leq c_1(n) [K(s):\Q] $, with $c_1(n)>0$ a constant depending only on $n$,
		
		\item $\deg(R_{s,\infty})\leq 8[L_s:\Q]$, and 
		
		\item the relation $R_{s,\infty}(\mathcal{Y}(x(s)))=0$ corresponding to the above polynomial is ``non-trivial''.	
\end{enumerate}	\end{prop}

\begin{remarks}We note that the points with unlikely isogenies for which all of the isogenous coordinates are singular are for all practical reasons dealt with by the work of Daw and Orr in \cite{daworr4}.
\end{remarks}

\begin{proof}Let $i_1\ld i_4$ be the four isogenous coordinates of $\CX_s$ and let us write $\CE_{i_j,0}$ for the fibers of the various connected N\'eron models at $s_0$. We assume without loss of generality that $i_1<i_2\leq i_3<i_4$.\\
	
	\textbf{Step 0: The issue of $v$-adic proximity} Here we may argue exactly as in the proof of Proposition $4.1$ of \cite{papaseffbrsieg}. We give a short summary of the ideas. If $s$ is $v$-adically close to $0$ with respect to some fixed archimedean place $v$, we get that $s$ is in some small archimedean disc around one of the $\xi_t$. Let $\lambda$ be the smallest element in the equivalence class of $t$, for the equivalence relation defined in \Cref{section:goodcovers}. Writing $s_t:=\sigma^{-1}_t(s)$, Lemma $5.4$ of \cite{daworr4} then gives an isogeny $\CX_{t,s_t}\sim \CX_{\lambda,s_{\lambda}}$. In particular $s_{\lambda}:=\sigma_\lambda^{-1}(s)$ will be a point with unlikely isogenies as those considered in \Cref{section:toycase}.\\
	
	\textbf{Step 1: The local factors} The rest of the proof goes as in Proposition $4.1$ of \cite{papaseffbrsieg}. Namely, to each $v$ as above we can assign one ``local factor'' $R_{s,v}$ for which $\iota_v(R_{s,v}(\mathcal{Y}_{\lambda}(x(s_{\lambda}))))=0$. Note here that $x(s_\lambda)=x(\sigma{^-1}_{\lambda}(s))=x(s)$ by the definition of the group $\aut_x(S)$, see \Cref{reduction2}. 
	
	The assumption that not all isogenous coordinates of $\CX_s$ are singular for $S'$ and the definition of $G_{ZP}$-admissibility shows that we are in either of the following situations:\\
	
	\textbf{Case 1:} $i_2=i_3$ and $\CE_{i_1,0}\times \CE_{i_2,0}\times \CE_{i_4,0}\simeq \G_m^2\times E$ with $E$ CM.\\
	
	The local factors $R_{s,v}$ in this case will be those constructed in \Cref{localfactorse2cm}.\\
	
	\textbf{Case 2:} $i_2=i_3$ and $\CE_{i_1,0}\times \CE_{i_2,0}\times \CE_{i_4,0}\simeq \G_m\times E\times E'$ with $E$, $E'$ both CM.\\
	
	The local factors $R_{s,v}$ are those constructed in \Cref{localfactors2ndway}.\\
	
	\textbf{Case 3:} $i_2\neq i_3$ and $\CE_{i_1,0}\times \CE_{i_2,0}\times \CE_{i_3,0}   \times \CE_{i_4,0}\simeq \G_m^3\times E$ with $E$ CM.\\
	
	The local factors $R_{s,v}$ in this case will be those constructed in \Cref{localfactors4}.\\
	
	\textbf{Case 4:} $i_2\neq i_3$ and $\CE_{i_1,0}\times \CE_{i_2,0}\times \CE_{i_3,0}   \times \CE_{i_4,0}\simeq \G_m^2\times E\times E'$ with $E$ and $E'$ both CM.\\
	
	There are two subcases here. If two of the isogenous coordinates, say $i_3$ and $i_4$, are CM then the local factors are those defined by \Cref{localfactors2ndway}.
	
	On the other hand, if none of the pairs of isogenous coordinates are both CM, we need to use the local factors $R_{s,v}$ of \Cref{localfactors4}.\\

	\textbf{Case 5:} $i_2\neq i_3$ and $\CE_{i_1,0}\times \CE_{i_2,0}\times \CE_{i_3,0}   \times \CE_{i_4,0}\simeq \G_m\times E\times E'\times E''$ with $E$, $E'$, and $E''$ all CM.\\
	
	In this case at least one of the pairs of isogenous coordinates are both CM. Thus, we can use the local factors $R_{s,v}$ of \Cref{localfactors2ndway}.\\
	
	\textbf{Case 6:} all of the coordinates $i_j$ are CM.\\
	
	The local factors are those defined by \Cref{localfactors2ndway}.\\
	
	\textbf{Step 3: The polynomial $R_{s,\infty}$} The definition of $R_{s,\infty}$ and the proof of its properties follow exactly as in the proof of Proposition $4.1$ of \cite{papaseffbrsieg}.
	
	In more detail, one defines \begin{equation}\label{eq:finalpoly}
		R_{s,\infty}(X^{(\lambda)}_{i,j,k}):=\prod_{\underset{s \text{ is } v\text{-adically close to }0}{v\in \Sigma_{L_s,\infty}}}^{}R_{{s,v}}(X^{(\lambda(v))}_{i,j,k}),
	\end{equation}note here that $\lambda$ will obviously depend on $v$. Parts $1$-$3$ then follow trivially by construction.
	
	For the non-triviality of $R_{s,\infty}$, we note that the fact that for all archimedean places $v$ we have $R_{s,v}\not\in I_0$ is enough, because $I_0$ is prime. Indeed, $V(I_0)$ is just $SL_2^{|\Lambda|\cdot k_0}$, for some $k_0\in \N$, in some affine space over $\bar{\Q}$, which is irreducible. It is thus enough to check that $I_0$ is a radical ideal which follows easily from Proposition $5.11$ of \cite{daworr5}.
\end{proof}

\begin{lemma}\label{localfactors4}
	Let $f':\CX'\rightarrow S'$ be a $G_{ZP}$-admissible semiabelian scheme with $n=4$. Let $s\in S(\bar{\Q})$ be some point with unlikely isogenies and let $L_s$ be its associated field of coefficients. Assume that $s$ is $v$-adically close to $\xi_\lambda$ with respect to some archimedean place $v\in \Sigma_{L_s,\infty}$ and that the following hold 
	\begin{enumerate}[label=(\roman*)]
		\item there are isogenies $\phi_1:\CE_{1,s}\rightarrow \CE_{2,s}$ and $\phi_2:\CE_{3,s}\rightarrow \CE_{4,s}$, and 
		
		\item either of the following holds
		\begin{enumerate}[label=(\alph*),ref=\theenumi(\alph*)]
			\item\label{lemma4case1} $1$ and $3$ are singular coordinates and the rest are CM, or
			
			\item\label{lemma4case2} $1$, $2$, $3$ are all singular coordinates while $4$ is a CM coordinate for $S'$.
		\end{enumerate}
	\end{enumerate}
	
	Then, there exists $R_{s,v}\in L_s[X^{(\lambda)}_{i,j,k}]$ such that the following hold\begin{enumerate}
		\item $\iota_v(R_{s,v}(\mathcal{Y}_{\lambda}(x(s))))=0$,
		
		\item $R_{s,v}$ is homogeneous of degree $\deg(R_{s,v})\leq 4$, and
		
		\item $R_{s,v}\notin I_0\leq L_s[X^{(\lambda)}_{i,j,k}]$, where $I_0$ is the ideal defined in \Cref{trivialrels}.
	\end{enumerate}
\end{lemma}
\begin{proof}Let us first assume that we are in \ref{lemma4case1}. 
	
	We work as in \Cref{section:toycase} in the ``first way'' of creating relations among the isogenous pair $\CE_{1,s}$ and $\CE_{2,s}$. We then get, as before $H_1^{(1)}$ and $H_2^{(1)}$ such that \eqref{eq:gm2xcm1} holds. In particular either $H_j^{(1)}=0$ for some $j$ or $H^{(1)}_1\cdot H^{(1)}_2=\frac{r_1s_1}{2\pi i}$ with $r_1s_1\neq0$.
	
	If $H_j^{(1)}=0$ we are done as in the proofs of earlier similar results. 
	
	Now working with the pair $\CE_{3,s}$ and $\CE_{4,s}$ we again get as before $H_1^{(2)}$ and $H_2^{(2)}$ such that \eqref{eq:gm2xcm1} holds with $r_2$, $s_2\in \Z$. Once again if $r_2=0$ or $s_2=0$ we are done as before. If on the other hand $r_2s_2\neq0$ we get $H^{(2)}_1\cdot H^{(2)}_2=\frac{r_2s_2}{2\pi i}$. 
	
	Assume from now own that $r_1s_1r_2s_2\neq 0$ so that we have\begin{equation}
		r_2s_2H^{(1)}_1\cdot H^{(1)}_2=r_1s_1H^{(2)}_1\cdot H^{(2)}_2.
	\end{equation}Then by similar arguments as in \Cref{localfactorse2cm} we get a polynomial $R_{s,v}$ that is homogeneous of degree $4$ and satisfies all of the properties that we want.\\
	
	Let us now assume that we are in \ref{lemma4case2}. By working with the isogenous pair $\CE_{3,s}$ and $ \CE_{s,4}$ we get on the one hand the same relations as in the previous case. Namely, reducing from above to the case $r_2s_2\neq 0$, we have \begin{equation}
		H^{(2)}_1\cdot H^{(2)}_2=\frac{r_2s_2}{2\pi i}.
	\end{equation}

Let us now work with the isogenous pair $\CE_{1,s}$ and $\CE_{2,s}$. Working as in the beginning of \Cref{section:gm2xe} and with the same notation for the various matrices as used there, we get that	\begin{equation}\label{eq:gmxgm1}
	\begin{pmatrix}a_1&0\\b_1&c_1\end{pmatrix} \begin{pmatrix}h^{(1)}_{i,j}	\end{pmatrix}\begin{pmatrix} d_1&e_{0,1}\\d'_1&e'_{0,1}	\end{pmatrix}=		
	\begin{pmatrix} h^{(2)}_{i,j}	\end{pmatrix}   \begin{pmatrix} d_2&e_{0,2}\\d'_2&e'_{0,2}	\end{pmatrix} \begin{pmatrix}p_1&q_1\\r_1&s_1	\end{pmatrix}.
\end{equation}

Arguing as in the ``first way'' of extracting relations described in \Cref{section:gm2xe} one again ends up with equations of the form \begin{equation}
	(H_1^{(1)},H_2^{(1)})=(\frac{r_1}{2\pi i}, \frac{s_1}{2\pi i}),
\end{equation}where $H_j^{(1)}$ are polynomials in the $h_{i,j}^{(k)}$ for $k=1$, $2$. These are nothing but a recreation of equation $(12)$ in \cite{daworr4}. 

We can then associate to $r_2s_2H^{(1)}_1\cdot H^{(1)}_2=r_1s_1H^{(2)}_1\cdot H^{(2)}_2$ a polynomial $R_{s,v}$ that will satisfy the conditions we want. The fact that only the first columns of the period matrices $\CP_{k,\iota_v}$ for $k=1$ and $2$ will appear follows from the construction of the $H_j^{(1)}$ as in the proof of \Cref{localfactorse2cm}.
\end{proof}

    
\section{Height bounds and Zilber-Pink}\label{section:heightbounds}

We conclude our exposition by establishing the height bounds needed for the large Galois orbits hypothesis and record some cases of the Zilber-Pink conjecture that follow from those, following work of Habegger-Pila \cite{habeggerpila1} and Daw-Orr \cite{daworr4}.

\subsection{Heights bounds}

We first need to rule out possible proximity of our points with unlikely isogenies to $s_0$ with respect to finite places. To do this we work as in \cite{papaseffbrsieg} assuming, which we may, that the technical Assumption $2.22$ there holds for our G-functions.

\begin{lemma}\label{padicproxzp}Let $f':\CX'\rightarrow S'$ be a $G_{ZP}$-admissible semiabelian scheme. Let $s\in S(\bar{\Q})$ be a point with unlikely isogenies and field of coefficients $L_s$.
	
	Assume that for one of the pairs of isogenous coordinates, say $i_1$ and $i_2$, of $\CX_s$ one of them is CM for $S'$ and the other one is singular for $S'$. Then if $v\in \Sigma_{L_s,f}$ is some finite place of $L_s$, the point $s$ is not $v$-adically close to $s_0$.
\end{lemma}
\begin{proof}
	By the same argument as in Lemma $5.1$ of \cite{papaseffbrsieg} we know that the special fiber of the connected N\'eron model of $\CX_{s}\times _{K(s)}L_{s,v}$ would be the same as that of $\CX_{0}\times _{K}L_{s,v}$.

	By Corollary $7.2$ of \cite{silvermanell} we also know that $\CE_{1,s}\times_{K(s)}L_{s,v}$ and $\CE_{2,s}\times_{K(s)} L_{s,v}$ will have the same type of reduction at $v$. By assumption we then have a contradiction since one of these will be $\G_{m,\kappa(v)}$, where $\kappa(v)$ here is the respective residue field, while the other one will be an elliptic curve over $\kappa(v)$. 
\end{proof}

\begin{theorem}\label{heightboundzp} Let $f':\CX'\rightarrow S'$ be a $G_{ZP}$-admissible semiabelian scheme. Then there exist effectively computable constants $c_1$ and $c_2$ such that for all $s\in S(\bar{\Q})$ that have unlikely isogenies and such that not all of the isogenous coordinates of $\CX_s$ are singular for $S'$, while all smooth isogenous coordinates of $\CX_s$ are CM for $S'$, we have \begin{equation}\label{eq:htbdao} h(s)\leq c_1[K(s):\Q]^{c_2}. \end{equation} \end{theorem}
\begin{proof}
	The proof is identical to that of Proposition $5.2$ in \cite{papaseffbrsieg}, replacing the usage of Proposition $4.1$ there by our \Cref{archrelastionsisogenies} and the usage of Lemma $5.1$ there by our \Cref{padicproxzp} respectively.
\end{proof}

\begin{remark}
	As in the discussion in \Cref{section:heightreduction} knowing \Cref{heightboundzp} for $G_{ZP}$ admissible curves is equivalent to knowing the height bound without any assumptions apart from:\begin{enumerate}
		\item the ``boundary-intersection'' conditions are the same as those in \Cref{heightboundzp}, and
		\item the image of the morphism $\iota:S\rightarrow Y(1)^n$, induced from the $n$ schemes over our curve $S$, is not contained in any proper special subvariety.
	\end{enumerate}
\end{remark}

\subsection{Some cases of the Zilber-Pink Conjecture}\label{section:applizp}

The strategy to reduce the Zilber-Pink conjecture for curves in $Y(1)^n$ to height bounds for isogenous points analogous to those that appear in \Cref{heightboundzp} already appears in \cite{daworr4}, based on work of Habegger and Pila in \cite{habeggerpila1}. 

Using the same arguments as in Proposition $5.12$ of \cite{daworr4} one can establish the following:
\begin{cor}[Large Galois Orbits for Zilber-Pink]\label{lgozp}
	Let $Z\subset Y(1)^n$ be an irreducible Hodge generic curve defined over $\bar{\Q}$ and let $K$ be  a field of definition of $Z$.
	
	Then there exist positive constants $c_3$, $c_4$ such that for every point $s\in Z(\bar{\Q})$ for which $\exists$ $\{i_1,i_2\}$, $\{i_3,i_4\}\subset \{1\ld n\}$ with $i_1\neq i_2$, $i_3\neq i_4$ and $\{i_1,i_2\}\neq\{i_3,i_4\}$ that are such that\begin{enumerate}
		\item $\exists$ $M$, $N$ with $\Phi_M(s_{i_1},s_{i_2})=\Phi_N(s_{i_3},s_{i_4})=0$, 
		
		\item $s_{i_1}$, $s_{i_3}$ are not singular moduli, and 
		
		\item $i_1,\dots,i_4$ are either singular or CM coordinates for $Z$,
	\end{enumerate}we have \begin{equation}
		[K(s):K]\geq c_3 \max\{M,N\}^{c_4}.
\end{equation}\end{cor}
\begin{proof}
	We simply note here the differences needed to adjust the proof of Proposition $5.12$ of \cite{daworr4} to our setting. 
	
	We adopt the notation of the proof of Theorem $1.5$ in \cite{papaseffbrsieg} finding a semiabelian scheme $f':\CX'\rightarrow S'$ that is $G_{ZP}$-admissible and such that $S'$ is a finite \'etale cover of $Z$.
	
	We can then apply \Cref{heightboundzp} to find $c_1$, $c_2$ with \begin{center}
		$h(x   (\tilde{s}))\leq  c_1[K(s): \Q]^{c_2}$
	\end{center} for all preimages $\tilde{s}$ in $S$ via $g$ of any such point $s\in Z(\bar{\Q})$.
	
	Letting $\rho_i$ be as in the previous proof, we recover the respective inequalities in the proof of Prop. $5.12$ in \cite{daworr4}, upon which stage we finish by using the isogeny estimates of Gaudron-R\'emond \cite{gaudronremond}.
\end{proof}

Given the above we can conclude from \cite{habeggerpila1} the following Zilber-Pink-type statement:
\begin{theorem}\label{zpgen}
	Let $C\subset Y(1)^n$ be an irreducible Hodge generic curve defined over $\bar{\Q}$. Let \begin{center}
		$J_1:=\{1\leq i\leq n: i\text{ is a singular coordinate for } C\}$ and $J_2:=\{1\leq i\leq n: i\text{ is a CM coordinate for } C\}$,
	\end{center}and set $J_C:=((J_1\cup J_2)\times (J_1\cup J_2))\backslash\{(i,i):1\leq i\leq n\}\subset\N^2$.
	Then the set 	
	\begin{center}
		$	\{ s\in C(\C):\exists N, M\text{ such that } \Phi_N(s_{i_1}, s_{i_2})= \Phi_M(s_{i_3}, s_{i_4})=0,(i_1,i_2)\neq (i_3,i_4),\text{ and }(i_1,i_2),(i_3,i_4)\in J_C \}$
	\end{center}is finite. 
\end{theorem}

	\bibliographystyle{alpha}
	\bibliography{biblio}

\end{document}